\documentclass[12pt]{article}

\usepackage{amsfonts}
\usepackage{amsmath}
\usepackage{amsthm} 
\usepackage{epsfig}
\usepackage{graphics}
\usepackage{pslatex}

\usepackage{latexsym}
\usepackage{amssymb,exscale}
\usepackage[all]{xy}
\usepackage{hyperref}
\usepackage{verbatim} 

\theoremstyle{plain} 
    \newtheorem{theorem}{Theorem}
    \newtheorem{lemma}[theorem]{Lemma}
    \newtheorem{proposition}[theorem]{Proposition}

\theoremstyle{definition} 

    \newtheorem{remark}[theorem]{Remark}
    \newtheorem{example}[theorem]{Example}

\def\newblock{\hskip .11em plus .33em minus .07em}

\def\eps{\epsilon}

\def\lam{\lambda}

\def\ra{\rightarrow}

\def\C{\mathbb{C}}

\def\GG{{\mathcal G}}

\def\R{\mathbb{R}}

\def\e{{\bf e}}

\def\summ{\sum\limits}
\def\intt{\int\limits}
\def\prodd{\prod\limits}

\def\tends{\rightarrow}

\def\l{\left}
\def\r{\right}
\def\<{\langle}
\def\>{\rangle}

\def\mb{\mbox}

\def\bar{\overline}

\def\given{\left.\vphantom{\hbox{\Large (}}\right|}

\newcommand\tr{{\mbox{\rm tr}}}

\newcommand\mnote[1]{} 
\newcommand\be{\begin{equation*}}

\newcommand\ee{\end{equation*}}

\newcommand\ben{\begin{equation}}
\newcommand\een{\end{equation}}
\newcommand\bes{\begin{eqnarray*}}
\newcommand\ees{\end{eqnarray*}}

\newcommand{\supp}{\mbox{\rm supp}}

\newcommand{\sm}{{\raise0.3ex\hbox{$\scriptstyle \setminus$}}}

\def\mb{\mbox}
\def\l{\left}
\def\r{\right}

\def\lam{\lambda}

\def\eps{\epsilon}
\def\tends{\rightarrow}

\renewcommand{\phi}{\varphi}

\def\CHI{\mathchoice%
{\raise2pt\hbox{$\chi$}}%
{\raise2pt\hbox{$\chi$}}%
{\raise1.3pt\hbox{$\scriptstyle\chi$}}%
{\raise0.8pt\hbox{$\scriptscriptstyle\chi$}}}
\def\smalloplus{\raise1pt\hbox{$\,\scriptstyle \oplus\;$}}

\author{Alice Guionnet\thanks{
UMPA, CNRS  UMR 5669, ENS Lyon, 46 all\'ee d'Italie,
69007 Lyon, France. {aguionne@umpa.ens-lyon.fr}.
This work was partially supported by the ANR project
ANR-08-BLAN-0311-01.}{ }, \;
 Manjunath Krishnapur\thanks{Department of Mathematics,
Indian Institute of Science, 
Bangalore - 560012, India.
{manju@math.iisc.ernet.in}}\;
{ }and Ofer Zeitouni\thanks{School of Mathematics,
University of Minnesota and Faculty of Mathematics,
Weizmann
Institute, POB 26, Rehovot 76100, Israel. {zeitouni@math.umn.edu}. 
The work of this author was partially
supported by NSF grant DMS-0804133 and by a grant from the Israel 
Science Foundation.}}
\title{The single ring theorem}

\begin{document}
 \bibliographystyle{abbrv}
\date{ {\small
September 10, 2009. Revised June 15 and October 16, 2010 and October 28, 2010}}
\maketitle


\noindent{\bf 
 Mathematics Subject of Classification :} {15A52 (46L50,46L54)}

\begin{abstract}
We study the empirical measure $L_{A_n}$
 of the eigenvalues of  non-normal square matrices
of the form $A_n=U_nT_nV_n$ with $U_n,V_n$ independent Haar distributed 
on the unitary group and $T_n$ real diagonal. We show that when the
 empirical measure of the eigenvalues of $T_n$ converges, and $T_n$
satisfies some technical conditions,  $L_{A_n}$
converges
towards a rotationally invariant measure $\mu$ on the complex plane whose
support is a single ring.
In particular, we provide a complete 
proof of Feinberg-Zee single ring theorem \cite{FZ}.  We also consider
the case where $U_n,V_n$ are independent Haar distributed 
on the orthogonal group.

\end{abstract}

\section{The problem}
Horn~\cite{horn} asked the question of describing the eigenvalues 
of a square matrix with prescribed singular values. 
If $A$ is a $n\times n$ matrix with singular values 
$s_{1}\ge \ldots \ge s_{n}\ge 0$ and eigenvalues 
$\lambda_{1},\ldots ,\lambda_{n}$ in decreasing order of absolute values, 
then the inequalities 
\begin{equation}
	\label{a-prioribound}
	\prodd_{j=1}^{k}|\lambda_{j}| \le \prodd_{j=1}^{k}s_{j}, \mb{ if }k< n 
 \qquad \mb{ and }\qquad \prodd_{j=1}^{n}|\lambda_{j}| = \prodd_{j=1}^{n}s_{j}
\end{equation}
were shown by Weyl~\cite{weyl} to hold. Horn established that these 
were all the relationships between singular values and eigenvalues. 

In this paper we study the natural probabilistic version of this problem and 
show that for ``typical matrices'', the singular values almost 
determine the eigenvalues. To frame the problem precisely, 
fix $s_1\ge \ldots \ge s_{n}\ge 0$ and consider $n\times n$ 
matrices with these singular values. They are of the form $A=PTQ$, 
where $T$ is diagonal with entries $s_{j}$ on the diagonal, 
and $P,Q$ are arbitrary unitary matrices. 

We make $A$ into a random matrix by choosing $P$ and $Q$ independently 
from Haar measure on $\mathcal{U}(n)$, the unitary group of 
$n\times n$ matrices, and independent from $T$. Let $\lambda_{1},\ldots ,\lambda_{n}$ be 
the (random) eigenvalues of $A$. The following natural questions arise.
\begin{enumerate}
\item Are there deterministic or random sets $\{s_{j}\}$, 
	for which one can find the exact distribution of $\{\lambda_{j}\}$? 
\item Let $L_{S}=\frac{1}{n}\sum_{j=1}^{n}\delta_{s_{j}}$ and 
	$L_{\Lambda}=\frac{1}{n}\sum_{j=1}^{n}\delta_{\lambda_{j}}$ 
	denote the empirical measures of $S=\{s_{j}\}$ 
	and $\Lambda=\{\lambda_{j}\}$. 
	Suppose $S_{n}$ are sets of size $n$ such that $L_{S_{n}}$ 
	converges weakly to a probability measure $\theta$ 
	supported on $\R_{+}$. Then, does $L_{\Lambda}$ 
	converge to a deterministic measure $\mu$ on the complex plane? 
	If so,  how is the measure $\mu$ determined by $\theta$?
\item For finite $n$, for fixed $S$, is $L_{\Lambda}$ {\em concentrated} 
	in the space of probability measures on the plane?
\end{enumerate} 
In this paper, we concentrate on the second question and
answer it in the affirmative, albeit 
with some restrictions. 
In this context,
we note that Fyodorov and Wei \cite[Theorem 2.1]{FW}
gave a formula for the mean  eigenvalues 
density of $A$, yet in terms of a large sum 
which does not offer an easy handle on asymptotic properties 
(see also \cite{FS} for the case where
$T$ is a projection). The authors of \cite{FW} explicitely state
the second
question as an open
problem.

Of course,
questions 1--3. above are not new, and have been
studied in various formulations. We now describe a partial
and necessarily 
brief history of what is known concerning questions 1. and 2.;
partial results concerning question 3. will be discussed elsewhere.

The most famous case of a positive answer to question 1. is the
{\em Ginibre ensemble}, see \cite{ginibre}, and its asymmetric variant,
see \cite{sommerslehman}. (There are some pitfalls in the standard
derivation of Ginibre's result. We refer to
 \cite{peresalsbook} for a discussion.)
Another situation is the truncation of random unitary matrices, described in
\cite{zysom}. 

Concerning question 2., the convergence of the empirical measure of eigenvalues
in the Ginibre ensemble (and other ensembles related
to question 1.) is easy to deduce from the explicit formula
for the joint distribution of eigenvalues. 
Generalizations of this convergence in the absence of
such explicit formula, for matrices with iid entries,
is covered under {\it Girko's circular law}, which is described  
in \cite{girko}; the circular law was proved under some conditions 
in \cite{bai} and finally, in full generality, in
\cite{gotzetikhomirov}
and \cite{taovu}.
Such matrices, however, do not possess the invariance properties
discussed in connection of question 2.
The {\em single ring theorem} of
Feinberg and Zee \cite{FZ} is, to our knowledge, the first 
example where a partial answer to this question is offered.
(Various issues of convergence are glossed over in \cite{FZ} and, as it
turns out, require a significant effort to overcome.) 
As we will see
in Section \ref{sec-aux},
the asymptotics of the spectral measure appearing in question
2. are described by the Brown
measure of $R$-diagonal operators. 
(The Brown measure is a continuous analogue of the spectral distribution of
non-normal operators, introduced in \cite{brown}.)  
$R$-diagonal operators were introduced by  Nica and Speicher \cite{NS}
in the context of free probability; they represent the weak*-limit 
(or more precisely, the limit in $*$-moments) of operators of the form 
$UT$ with $U$ unitary with size going to infinity 
and $T$ diagonal,  
and were
intensively studied in the last decade within the theory of free
probability, in particular in connection with the problem of classifying
invariant subspaces \cite{HS, HS2}.

\section{Limiting spectral density of a non-normal matrix}
Throughout, for a probability measure $\mu$ supported
on $\R$ or on $\C$, we write
$G_\mu$ for its Stieltjes transform, that is
$$G_\mu(z)=\int \frac{\mu(dx)}{z-x}\,.$$
$G_\mu$ is analytic off the support of $\mu$.
We let ${\mathcal H}_n$ denote the Haar measure 
on the $n$-dimensional unitary group $\mathcal{U}(n)$.
Let $\{P_n,Q_n\}_{n\geq 1}$ denote a sequence of independent,
${\mathcal H}_n$-distributed matrices. Let $T_n$ denote a sequence
of diagonal matrices, independent of $(P_n,Q_n)$,
 with real positive entries 
$S_{n}=\{s^{(n)}_{i}\}$ 
on the diagonal, and introduce the {\em empirical measure}
of the {\em symmetrized} version of $T_n$ as 
$$L_{S_{n}}=\frac1{2n} \sum_{i=1}^n [\delta_{s^{(n)}_i}+
\delta_{-s^{(n)}_i}]\,.$$ 
We write $G_{T_n}$ for $G_{L_{S_n}}$.
For a measure $\mu$ supported on $\R_+$, we write 
$\tilde \mu$ for its {\em symmetrized version}, that is,
for any $0<a<b<\infty$,
$$\tilde \mu([-a,-b])=\tilde \mu([a,b])=\frac12 \mu([a,b])\,.$$

Let $A_n=P_n T_nQ_n$,
let
$\Lambda_{n}=\{\lambda_{i}^{(n)}\}$ denote the set of eigenvalues of $A_n$,
and set 
$$L_{A_{n}}=\frac1n \sum_{i=1}^n \delta_{\lambda^{(n)}_i}\,.$$ 
We refer to $L_{A_n}$ as the empirical spectral distribution (ESD) of
$A_n$.
(Note that the law of
$L_{A_n}$ does not change if one considers $P_n T_n$ instead
of $P_nT_nQ_n$, since if $P_nT_nQ_n w=\lambda w$
for some $(w,\lambda)$ then, with $\bar P_n=Q_n P_n$ and $v=Q_n w$, 
it holds that 
$\bar P_n T_n v=\lambda v$, and $\bar P_n$ is again Haar distributed.) 
Finally,
for any matrix $A$, we set $\|A\|$ to denote the $\ell^2$
operator-norm of $A$, that is, its largest singular value.

To state our results, we recall the notion of {\it free convolution}
of probability measures on $\R$, introduced by Voiculsecu. 
For a compactly supported 
probability measure on $\mu$, define the 
formal power series
$G_\mu(z)=\sum_{n\geq 0} \int x^n d\mu(x) z^{-(n+1)}$, and let
$K_\mu(z)$ denote its inverse in a neighborhood of infinity, 
satisfying $G_\mu(K_\mu(z))=z$.
The {\it R-transform} of $\mu$ is the the function $R_\mu(z)=K_\mu(z)-1/z$.
The moments of $\mu$ (and therefore $\mu$ itself, since it is compactly
supported) can be recovered from the knowledge of $K_\mu$, and therefore
from $R_\mu$, by
a formal inversion of power series. For a pair of compactly supported 
probability measures $\mu_1,\mu_2$, introduce the {\it free convolution}
$\mu_1\boxplus \mu_2$ as the (compactly supported) probability measure
whose R-transform is $R_{\mu_1}(z)+R_{\mu_2}(z)$. (That this defines
indeed a probability measure needs a proof, see \cite[Section 5.3]{AGZ}
for details and background.)

For $a\in \R_+$, introduce the symmetric Bernoulli measure
$\lambda_a=\frac12(\delta_a+\delta_{-a})$ with atoms at $\{-a,a\}$.
All our main results, Theorem \ref{cor-FZ} and Propositions
\ref{cor-D} and \ref{corplus},
will be derived from the following technical
result.
\begin{theorem}
	\label{main-theo}
	Assume $\{L_{T_n}\}_n$ converges weakly
	to a probability measure $\Theta$ compactly supported
	on $\R_+$.
	Assume further 
	\begin{enumerate}
		\item
			There exists a constant $M>0$
	so that 
	\begin{equation} \label{bound}\lim_{n\to\infty} P(\|T_n\|>M)=0\,.
	\end{equation}
\item There exist a sequence of events $\{{\cal G}_n\}$ with
	$P({\cal G}_n^c)\to 0$ and
	constants $\delta,\delta'>0$ so that for Lebesgue almost any $z\in \C$,
	with $\sigma_n^z$ the minimal singular value
	of $zI-A_n$,
	\begin{equation}
		\label{eq-300609b}
		E({\bf 1}_{{\cal G}_n}
		{\bf 1}_{\{\sigma_n^z<n^{-\delta}\}}
		(\log \sigma_n^z)^2)<\delta'\,.
	\end{equation}
\item There exist  constants $\kappa,\kappa_1>0$ such that
	\begin{equation}
		\label{eq-denbound}
		|\Im G_{T_n}(z)|\leq \kappa_1 \quad\mbox{\rm on}\quad 
		\{z:\Im(z)>n^{-\kappa}\}\,.
	\end{equation}
	\end{enumerate}
	Then the following hold.
\begin{itemize}
\item[a.] $L_{A_n}$ converges in probability to a limiting 
	probability measure
	$\mu_A$.
\item[b.] 
The measure $\mu_A$ possesses a radially-symmetric
density $\rho_A$ with respect 
to the Lebesgue measure on $\C$, satisfying
$\rho_A(z)=\frac{1}{2\pi}\Delta_z (\int \log |x| d\nu^z(x))$, 
 where
$\Delta_z$ denotes the Laplacian with respect to the variable $z$ 
and $\nu^z
:=\tilde \Theta\boxplus \lambda_{|z|}$\;.
\item[c.] The support of $\mu_A$ 
	is a single ring: there exist constants
	$0\leq a<b<\infty$ so that
	$$\supp \mu_A=\{re^{i\theta}: a\leq r\leq b\}\,.$$
	Further, $a=0$ if and only if $\int x^{-2} d\Theta(x)=\infty$.
\end{itemize}
\end{theorem}
See Remark \ref{remark-law} for an explicit characterization
of the free convolution appearing in Theorem \ref{main-theo},
and \cite[Ch. 5]{AGZ} for general background.
A different characterization of $\rho_A$, borrowed from
\cite{haageruplarsen} and instrumental
in the proof of part (c) of Theorem \ref{main-theo}, 
is provided in Remark \ref{rem-stransf} in Section
\ref{subsec-lim}.

\begin{remark}
	We do not believe that the conditions in Theorem
	\ref{main-theo} are sharp. In particular, we do not know whether
	condition 4, which prevents the existence of an atom in the support
	of $\tilde \Theta$, can be dispensed of;
	the example $T_n=I$ shows that it is certainly not necessary.
\end{remark}
\noindent
Theorem \ref{main-theo} is generalized to the case where $U_n,V_n$
follow
the Haar measure on the orthogonal group in Theorem
\ref{theo-orth}.
Note that, since for Lebesgue almost every $x\in \R$, the imaginary part 
of the Stieltjes transform of an absolutely continuous probability measure 
converges, as $z\to x+i\eps$,
towards the density of this 
measure at $x$, \eqref{eq-denbound} is verified as soon as 
$\tilde\Theta$ has a bounded continuous density.

As a corollary of 
Theorem \ref{main-theo}, we prove the Feinberg-Zee ``single ring theorem''.
\begin{theorem}
	\label{cor-FZ}
	Let $V$ denote a polynomial with positive leading 
coefficient. Let 
the $n$-by-$n$ complex matrix $X_n$ be distributed according to the law
$$ \frac{1}{Z_n} \exp(-n\tr \ V(XX^*))dX\,,$$
where $Z_n$ is a normalization constant and $dX$ the Lebesgue measure
on $n$-by-$n$ complex matrices. 
Let $L_{X_n}$ be the ESD
of
$X_n$. Then $\{L_{X_n}\}_n$ satisfies the conclusions of
Theorem \ref{main-theo} with $\Theta$ the unique minimizer of the
functional 
$${\cal J}(\mu):=\int V(x^2)d\mu(x)-\int\int\log|x^2-y^2|d\mu(x)d\mu(x)$$
on the set of probability measures  on $\mathbb R^+$. 
\end{theorem}
Theorem \ref{cor-FZ} will follow by checking that the assumptions
of Theorem \ref{main-theo} are satisfied for the spectral decomposition
$X_n=U_n T_n V_n$, see Section \ref{sec-corproof}.

The second hypothesis in Theorem \ref{main-theo}
may seem difficult to verify in general; we show in the next proposition 
that adding a small Gaussian matrix guarantees it.
\begin{proposition}
	\label{cor-D}
Let $(T_n)_{n\ge 0}$ be a sequence of matrices satisfying
the assumptions of Theorem \ref{main-theo}
except for  \eqref{eq-300609b} and assume that $\|T_n^{-1}\|$
is uniformly bounded. Let $N_n$ be  a $n\times n$
	matrix with independent (complex) Gaussian 
	entries of zero mean and covariance equal identity.
Let $U_n,V_n$ follow the Haar measure on unitary $n\times n$ matrices,
independently of $T_n,N_n$.
Then, the empirical measure of the eigenvalues 
of $Y_n:= U_nT_n V_n+ n^{-\gamma} N_n$ converges weakly in probability
to $\mu_A$ as in Theorem \ref{main-theo}  for any $\gamma\in (\frac
12,\infty)
$. 
\end{proposition}
\begin{example}\label{exa}
An example of sequence  $(T_n)_{n\ge 0}$ satisfying 
the hypotheses of Proposition \ref{cor-D} is given as follows: take
$\mu$  a  compactly supported 
 probability measure on $\mathbb R^{+}$.
Assume  the inverse $F^{-1}$ of the
distribution function $F(x)=\mu([0,x])$ is H\"older continuous 
and that  the imaginary part of the
Stieltjes transform of $\mu$  is  uniformly bounded
 on $\mathbb C^+$. 
Then the diagonal matrix $T_n$ with entries  
$$s_i^n=\inf\{ s:\mu([0,s])\ge \frac{i}{n}\},\quad 1\le i\le n\,,$$
satisfies the hypotheses of Proposition \ref{cor-D}.
\end{example}
A rather straightforward generalization of Theorem \ref{main-theo}
concerns the limiting spectral measure of $P_n+B_n$, where $P_n$ is
${\mathcal H}_n$ distributed and
the sequence of
$n\times n$ matrices $B_n$ converges in  $*$-moments to an operator
$b$ in a non-commutative probability space $({\mathcal A},\tau)$.
(The latter means 
 that for all polynomial $P$ in two non-commutative variables,
$$\lim_{n\ra\infty}\frac{1}{n}\tr\left(P(B_n,B_n^*)\right)=
\tau(P(b,b^*))\,,$$
which is the case if e.g
$B_n$ is self-adjoint, with spectral measure converging
 to a probability measure
$\Theta$, which is the law of a self-adjoint operator $b$.)
In particular,
for any $w\in \mathbb C$, the spectral measure 
of $T_n(w)=|wI-B_n|=\sqrt{(wI-B_n)(wI-B_n)^*}$ converges
to the law $\Theta_w$  of $|wI-b|$. 
By Voiculescu's theorem \cite[Theorem 3.8]{Vo91}, if the operator
norm of $B_n$ is uniformly bounded, then
the couple  $(B_n,P_n)$ converges in  $*$-moments 
towards $(b,u)$, a couple of operators living in 
a non-commutative probability space $({\mathcal A},\tau)$
which are free, $u$ being unitary.  The Brown measure 
$\mu_{b+u}$ is studied in \cite[Section 4]{bl01}.
\begin{proposition}\label{corplus}
Assume that $T_n(0)$ satisfies \eqref{bound} and that there exists a set 
$\Omega\subset \mathbb C$
with full Lebesgue measure so that  for all $w\in\Omega$, 
$T_n(w)$ satisfies 
\eqref{eq-denbound}. 
Let $N_n$ be  a $n\times n$
	matrix with independent (complex) Gaussian 
	entries of zero mean and covariance equal identity.
	Then, for any $\gamma>\frac{1}{2}$,
the spectral measure of $B_n +n^{-\gamma} N_n+P_n$ converges 
in probability to the Brown measure $\mu_{b+u}$  of $b+u$.
\end{proposition}
An example of matrices $B_n$ which satisfy the hypotheses of Proposition
\ref{corplus}
is given by the diagonal matrices $B_n=\mbox{diag}(s_1^n,\ldots, s_n^n)$
 with entries $s_i^n$ 
 satisfying the hypotheses of
 Example \ref{exa}. This is easily verified from 
the fact that the eigenvalues of $D_n(w)$ are given by
$(|w-s_1^n|,\ldots, |w-s_n^n|)$.

\subsection{Background and description of the proof}
The main difficulty in studying the ESD 
$L_{A_{n}}$ is that $A_n$ is not a normal matrix, that is 
$A_nA^{*}_n\not= A^{*}_nA_n$, almost surely. 
For normal matrices, the limit of 
ESDs can be found by the method of moments or by the method of 
Stieltjes' transforms. For non-normal matrices, the only known method of
proof is more indirect and follows an idea of Girko~\cite{girko} that 
we describe now (the details are  a little different from what is 
presented in Girko~\cite{girko} or Bai~\cite{bai}).

From Green's formula, for any polynomial 
$P(z)=\prod_{j=1}^{n}(z-\lambda_{j})$, we have
$$
 \frac{1}{2\pi}\int \Delta \psi(z) \log|P(z)| dm(z) = \summ_{j=1}^{n} 
 \psi(\lambda_{j}), \qquad \mb{ for any }\psi \in C_{c}^{2}(\C)\,,
$$
where $m(\cdot)$ denotes the Lebesgue measure on $\C$.
Applied to the characteristic polynomial of $A_n$, this gives
\begin{eqnarray*}
 \int \psi(z) dL_{A_{n}}(z)  &=& 
 \frac{1}{2\pi n} \intt_{\C} \Delta \psi(z) \log|\det(zI-A_n)| dm(z) \\
  &=& \frac{1}{4\pi n} \intt_{\C} \Delta \psi(z) \log \det(zI-A_n)(zI-A_n)^{*}
  dm(z)\,.
\end{eqnarray*}
It will be convenient for us to introduce the $2n\times 2n$ matrix
\begin{equation}
	\label{eq-hnz}	
 H_{n}^{z}:=\l[\begin{array}{cc} 0 & zI-A_n \\ (zI-A_n)^* & 0 \end{array}  \r].
\end{equation}
It may be checked easily that eigenvalues of $H_{n}^{z}$ are the positive and 
negative of the singular values of $zI-A_n$. 
Therefore, if we let $\nu_{n}^{z}$ denote the ESD of $H_{n}^{z}$, 
$$\int \frac{1}{y-x} d\nu_n^z(x)=\frac{1}{2n}\tr\left(
(y-H_{n}^{z})^{-1}
\right)\, ,$$ 
then 
$$
 \frac{1}{n}\log \det(zI-A_n)(zI-A_n)^{*} = 
 \frac{1}{n}\log \det |H_{n}^{z}| =2 \int_{\R} \log |x| d\nu_{n}^{z}(x)\,.
$$  
Thus we arrive at the formula
\begin{equation}\label{eq:girko}
\int \psi(z) dL_{A_{n}}(z) =  \frac{1}{2\pi} \intt_{\C} \Delta \psi(z) 
\int_{\R}\log |x| d\nu_n^z(x)d m(z)\,.\end{equation}
This is Girko's formula in a different form and its utility 
lies in the following attack on finding the limit of $L_{A_{n}}$.
\begin{enumerate}
\item Show that for (Lebesgue almost)
	every $z\in \C$, the measures $\nu_{n}^{z}$ 
	converge 
 weakly in probability 
to a measure $\nu^{z}$ as $n\tends \infty$,
	and identify the limit. 
	Since $H_{n}^{z}$ are Hermitian matrices, 
	there is hope of doing this by Hermitian techniques.
\item Justify that $\int \log|x| d\nu_{n}^{z}(x) \rightarrow \int \log|x| 
	d\nu^{z}(x)$ for (almost every) 
	$z$. But for the fact that ``$\log$'' 
	is not a bounded function, this would have followed from 
	the weak convergence of $\nu_{n}^{z}$ to $\nu^{z}$. 
	As it stands, this is the hardest technical part of the proof.
\item A standard weak convergence argument 
	is then used in order to convert the convergence 
	for (almost every) $z$ of $\nu_n^z$
	to a convergence of integrals over $z$. Indeed,
	setting $h(z):= \int \log|x| d\nu^{z}(x)$,
	we will get from \eqref{eq:girko} that 
	\begin{equation}
		\label{eq-h}
 \int \psi(z) dL_{A_{n}}(z)  \rightarrow \frac{1}{2\pi} 
 \intt_{\C} \Delta \psi (z) \; h(z) dm(z)\,. 
\end{equation}
\item Show that $h$ is smooth enough so that one can 
	integrate the previous equation by parts to get
\begin{equation}
\label{eq-290709a}
 \int \psi(z) dL_{A_{n}}(z)  \rightarrow \frac{1}{2\pi} 
 \intt_{\C}  \psi (z) \; \Delta h(z) dm(z)\,,  
\end{equation}
which identifies $\Delta h(z)$ as the density 
(with respect to Lebesgue measure) of the limit of $L_{A_{n}}$. 
\item Identify the function $h$ sufficiently precisely to be 
	able to deduce properties of $\Delta h(z)$. 
	In particular, show the {\bf single ring phenomenon}, 
	which states that the support of the limiting spectral measure 
	is a single annulus (the surprising part being that it 
	cannot consist of several disjoint annuli). 
\end{enumerate}
Girko's equation~\eqref{eq:girko} and these five steps give a 
general recipe for finding limiting spectral measures of 
non-normal random matrices. Whether one can overcome the technical 
difficulties depends on the model of random matrix one chooses. 
For the model of random matrices with i.i.d. entries having zero 
mean and finite variance, this has been achieved in stages by 
Bai~\cite{bai}, G\"{o}tze and Tikhomirov~\cite{gotzetikhomirov}, 
Pan and Zhou \cite{panzhou} and 
Tao and Vu~\cite{taovu}. While we heavily
borrow from that sequence, a major difficulty in the problem
considered here is that there is 
no independence between entries of the matrix
$A_n$. 
Instead, we will rely on properties 
of the Haar measure, and in particular on considerations
borrowed from free probability and the so 
called {\em Schwinger--Dyson} (or {\em master-loop}) equations.
Such equations were already the key  
to obtaining
fine estimates on the Stieltjes transform of Gaussian
generalized band matrices in \cite{HT}.  In
\cite{CGM}, they were used  to study the
asymptotics 
of matrix models on the unitary group. Our approach 
combines ideas of \cite{HT} to estimate Stieltjes transforms
and the necessary adaptations to unitary matrices as developped in \cite{CGM}.
The main observation is that one can reduce attention to the study 
of the ESD of matrices of the form $(T+U)(T+U)^*$ where $T$ is real diagonal 
 and $U$ is Haar distributed. In the limit (i.e., when $T$ 
and $U$ are replaced by operators in a $C^*$-algebra
that are freely independent,
with $T$ bounded and self adjoint and $U$ unitary), the limit
ESD has been identified by Haagerup and Larsen
\cite{haageruplarsen}. The Schwinger--Dyson equations give both
a characterization of the limit and, more important to us, a discrete
approximation that can be used to estimate the discrepancy between
the pre-limit ESD and its limit. These estimates play a crucial role
in integrating the singularity of the log in Step two above, but only once
an a-priori (polynomial) estimate on the minimal singular value has 
been obtained. The latter is deduced from assumption \ref{eq-300609b}.
In the context of the Feinberg--Zee single ring theorem, the latter
assumption holds due to an adaptation of 
the analysis of \cite{SST}.

\subsection*{Notation}
We describe our convention concerning constants. Throughout, by the word 
{\em constant} we mean quantities
that are independent of $n$
(or of the complex variables $z$, $z_1$). Generic 
constants
denoted by the letters $C$,$c$ or $R$, have
values that may change
from line to line, and they may depend on other parameters. Constants
denoted by $C_i$, $K$,   $\kappa$ and $\kappa'$
are fixed and do not change from line to line.
\section{An auxiliary problem: evaluation of $\nu^z$ and convergence rates}
\label{sec-aux}
Recall from the proof sketch described above that we are interested
in evaluating the limit $\nu^z$ of the ESD $L_n^z$ of the matrix
$H_n^z$, see \eqref{eq-hnz}. Note that 
$L_n^z$ is also the ESD of the matrix $\tilde H_{n}^{z}$ given by 
\begin{eqnarray}
	\label{eq-hnz1}	
\tilde H_{n}^{z}&:=&
  \l[\begin{array}{cc} 0 & Q_n \\ P_n^* & 0 
 \end{array}  \r]
 H_n^z
  \l[\begin{array}{cc} 0 & P_n \\ Q_n^* & 0 
 \end{array}  \r]\\
 &=&
\l[\begin{array}{cc} 0 & |z|W_n^z-T_n \\ (|z|W_n^z-T_n)^* & 0 
 \end{array}  \r]
\nonumber
 \,,
\end{eqnarray}
where $W_n^z=\bar z Q_n P_n/|z|$ is unitary and $\mathcal{H}_n$
distributed. Throughout, we will write $\rho=|z|$. We also will assume
in this
section that
the sequence $T_n$ is 
deterministic. We  are thus led to the study of the ESD for a sequence 
of matrices of the 
form
\begin{equation}
	\label{eq-070709a}
	{\bf Y}_n=\left(\begin{array}{cc}
0&B_n\\
B^*_n&0\\
\end{array}\right)
\end{equation}
with $B_n=\rho U_n+T_n$, $T_n$ being a real, diagonal
matrix of uniformly bounded norm, 
and $U_n$ a ${\mathcal H}_n$ unitary matrix.  
Because $\|T_n\|$ is uniformly bounded, it will be enough to consider 
throughout $\rho$ uniformly bounded.

We denote in short
\begin{equation}
	\label{eq-070609b}
	{\bf U}_n=\left(\begin{array}{cc}
0&U_n\\
0&0
\end{array}\right)\,,\quad {\bf U}^*_n=\left(\begin{array}{cc}
0&0\\
U^*_n&0
\end{array}\right)\,,\quad {\bf T}_n=\left(\begin{array}{cc}
0&T_n\\
T_n&0
\end{array}\right)\,.
\end{equation}

\subsection{Limit equations}
\label{subsec-lim}
We begin by deriving the limiting Schwinger--Dyson equations for the ESD 
of 
${\bf Y}_n$. Throughout this subsection, we consider a non-commutative
probability space 
$({\cal A}, *,\mu)$ on which a variable  $U$ lives
and 
where $\mu$ is a tracial state 
 satisfying the relations $\mu((UU^*-1)^2)=0$, $\mu(U^a)=0$ 
for $a\in \mathbb{Z}\setminus \{0\}$. 
 In the sequel, $1$ will
denote the identity in ${\cal A}$.
We refer to \cite[Section 5.2]{AGZ}
for definitions.

Let $T$ be a self-adjoint (bounded) element in ${\cal A}$,
with $T$ freely independent  with $U$. Recall the non-commutative derivative
$\partial$, defined on elements of $\mathbb{C}\langle T,U,U^*\rangle$ as 
satisfying 
the Leibniz rules
\begin{equation}\label{leibniz}
\partial (PQ)=\partial P\times (1\otimes Q)+(P\otimes 1)\times \partial Q\,,
\end{equation}
$$\partial U=U\otimes 1,\, \partial U^*=-1\otimes U^*,\, \partial T=0\otimes 0
\,.$$ 
(Here, $\otimes$ denotes the tensor product and 
we write $(A\otimes B) \times (C\otimes D)=(AC)\otimes(BD)$.)
$\partial$ is defined so that for any $B\in {\cal A}$
satisfying
$B^*=-B$, any $P\in \mathbb{C}\langle U,U^*, T\rangle$,
\begin{equation}
	\label{eq-new0928}P(Ue^{\epsilon B},  e^{-\epsilon B}U^*, T)=P(U,U^*,T)+\epsilon
\partial P (U,U^*,T)\sharp B+o(\epsilon)
\,,
\end{equation}
where we used the notation $A\otimes B\sharp C=ACB$.

By the invariance of $\mu$ under unitary conjugation, 
see \cite[Proposition 5.17]{VoiculescuAdvances}
or \cite[(5.4.31)]{AGZ}, we have the Schwinger--Dyson equation
\begin{equation}
	\label{eq-060709a}
	\mu\otimes \mu(\partial P)=0\,.
\end{equation}

 We continue to use the notation ${\bf Y}$, 
${\bf U}, {\bf U^*}$ and
${\bf T}$ in a way similar to
\eqref{eq-070709a} and \eqref{eq-070609b}. So,
we let ${\bf Y}=\rho({\bf U}+{\bf U}^*)+{\bf T}$
with
\begin{equation}
	\label{eq-070609bb}
	{\bf U}=\left(\begin{array}{cc}
0&U\\
0&0
\end{array}\right)\,,\quad {\bf U}^*=\left(\begin{array}{cc}
0&0\\
U^*&0
\end{array}\right)\,,\quad {\bf T}=\left(\begin{array}{cc}
0&T\\
T&0
\end{array}\right)\,.
\end{equation} 
We extend $\mu$ to the algebra generated by ${\bf U},{\bf U}^*$ and
${\bf T}$ by putting for any $A,B,C,D\in{\cal A}$,

$$\mu\left(\left(\begin{array}{cc}
A&B\\
C&D
\end{array}\right)\right):=\frac{1}{2}\mu(A)+\frac{1}{2}\mu(D)\, .$$
Observe that this extension is still tracial. 

The non-commutative derivative $\partial$ in 
\eqref{eq-new0928} extends naturally
to the algebra generated by the matrix-valued ${\bf U}, {\bf U}^*,
{\bf T}$, using the Leibniz rule \eqref{leibniz} together with the relations
\begin{equation}\label{diffp}
\partial{\bf U}={\bf U}\otimes p\,,
\quad \partial{\bf U^*}=- p\otimes {\bf U^*}
\,,\quad\partial{\bf T}=0\otimes 0\,,\end{equation}
where we denoted $p=\left(\begin{array}{ll}
0&0\\
0&1\\
\end{array}\right)$. In the sequel we will 
apply $\partial$ to analytic functions of ${\bf U}+{\bf U^*}$ and ${\bf T}$
 such as products
of Stieltjes functionals  of the form $\left(z-b{\bf U}-b{\bf U^*}-
a {\bf T}\right)^{-1}$ with $z\in {\mathbb C}\backslash {\mathbb R}$ and $a,b
\in \mathbb R$.
Such an extension is straightforward; $\partial$ continues to satisfy the 
Leibniz rule and, using the resolvent identity 
\begin{eqnarray*}
	&&\partial \left(z-b{\bf U}-b{\bf U^*}-a {\bf T}\right)^{-1}=
\\&&
b\left(z-b{\bf U}-b{\bf U^*}-a {\bf T}\right)^{-1} \left( {\bf U}\otimes p-
 p\otimes {\bf U^*}\right)\left(z-b{\bf U}-b{\bf U^*}-a {\bf T}\right)^{-1}
 \,,\end{eqnarray*}
 where
 $A(B\otimes C)D=(AB)\otimes(CD)$.
 Further, \eqref{eq-060709a} extends also in this context. 

Introduce the notation, for $z_1,z_2\in {\mathbb C}^+$, 
\begin{eqnarray}
	\label{eq-070809a}
G(z_1,z_2)&=&\mu\left( (z_1-{\bf Y})^{-1}(z_2-{\bf T})^{-1}\right)\nonumber\,,\\
G_U(z_1,z_2)&=&\mu\left( {\bf U}(z_1-{\bf Y})^{-1}(z_2- {\bf
T})^{-1}\right)\nonumber\,,\\
G_U(z_1)&=&\mu\left( {\bf U}(z_1-{\bf Y})^{-1}\right)\nonumber\,,\\
G_{U^*}(z_1,z_2)&=&\mu\left( {\bf U^*}(z_1-{\bf Y}
)^{-1}(z_2-{\bf T})^{-1}\right)\,,\\
G_{T}(z_1,z_2)&=&\mu\left( {\bf T}(z_1-{\bf Y}
)^{-1}(z_2 -{\bf T})^{-1}\right)\nonumber\,,\\
G(z_1)&=&\mu\left( (z_1-{\bf Y})^{-1}\right)\nonumber\,,\\
G_T(z_2)&=&\mu\left( (z_2 -{\bf T})^{-1}\right)\,.\nonumber
\end{eqnarray}
We apply the derivative $\partial$ to  the analytic function
$P=(z_1-{\bf Y})^{-1}(z_2- {\bf T})^{-1}{\bf U}\,,$
while noticing that, by \eqref{leibniz} and \eqref{diffp},
\begin{equation}\label{caf}
\partial P=P\otimes p+\rho (z_1-{\bf Y})^{-1} {\bf U}\otimes p P
-\rho (z_1-{\bf Y})^{-1}p \otimes {\bf U}^* P. 
\end{equation}
Applying \eqref{eq-060709a}, with $\mu(P)=G_U(z_1,z_2)$ and $\mu(p)=1/2$, we
find
\begin{equation}\label{cafcaf}
\frac{1}{2} G_U(z_1,z_2)=\rho \mu\left( (z_1-{\bf Y})^{-1}p\right)\mu({\bf U}^* P)
-\rho \mu\left((z_1-{\bf Y})^{-1} {\bf U}\right)\mu (pP)\,.
\end{equation}
Note that
$ Pp=P$ and thus  $\mu(pP)=\mu(P)$. Further,
for any  smooth function $Q$,
$\mu({\bf U}^* Q {\bf U})$ equals $ \mu( (1-p) Q)$
due to the traciality of $\mu$ and ${\bf U}{\bf U}^*=1-p$.
By symmetry (note that 
 $(1-p)(z_1-{\bf Y})^{-1}(z_2-{\bf T})^{-1}$ and $p(z_1-{\bf
Y})^{-1}(z_2-{\bf T})^{-1} $ are given by the same formula up
to replacing $({\bf U},{\bf U}^*)$ by $({\bf
U}^*,{\bf U})$, which has the same law) we get that
$\mu({\bf U}^*P)$ equals 
\begin{equation}
	\label{eq-new0928a}
	\mu((1-p)(z_1-{\bf Y})^{-1}(z_2-{\bf T})^{-1})
=\frac{1}{2}\mu((z_1-{\bf Y})^{-1}(z_2-{\bf T})^{-1})=
\frac12 G(z_1,z_2)\,.
\end{equation}
The first equality holds without the 
last factor $(z_2-{\bf T})^{-1}$, thus implying that
$\mu( (z_1-{\bf Y})^{-1} p)=\mu( (z_1-{\bf Y})^{-1})/2=G(z_1)/2$
and so we get
from \eqref{cafcaf}
that
\begin{equation}\label{eqU}
\frac{1}{2} G_U(z_1,z_2)=\frac{\rho}{4} G(z_1,z_2)G(z_1)
-\rho G_{U}(z_1,z_2)G_U(z_1)\,.
\end{equation} 
Noticing that $G_U(z_1)$ is the limit 
of $ z_2 G_U(z_1,z_2)$ as $z_2\to\infty$, 
 we find by \eqref{eqU} that
$$\frac{1}{2}G_U(z_1)=-\rho G_{U}(z_1)^2+\frac{\rho}{4} G(z_1)^2\,,$$
and therefore, as $G_U(z_1)$ goes to zero as $z_1\to\infty$, 
\begin{equation}\label{eq3}
G_U(z_1)= \frac{1}{2\rho}(-\frac{1}{2}+\sqrt{\frac{1}{4}
+ \rho^2G(z_1)^2})=\frac{1}{4\rho}(-1+\sqrt{1+4\rho^2G(z_1)^2})\,.
\end{equation}
Here, the choice of the branch of the square root is determined
by the expansion of $G_U(z)$ at infinity and the fact
that both $G(z)$ and $G_U(z)$ are analytic in $\mathbb C^+$.
This equation is then true for all $z_1\in\C^+$.

Moreover, by \eqref{eqU} and \eqref{eq3}, we get
\begin{equation}\label{eq33}
G_U(z_1,z_2)=\frac{\rho}{2}\frac{G(z_1,z_2)G(z_1)}{1+2\rho G_U(z_1)}
=\frac{\rho G(z_1,z_2)G(z_1)}{1+\sqrt{1+4\rho^2  G(z_1)^2}}\,.
\end{equation}
(Again, here and in the rest of this subsection, the proper branch of 
the square root is determined by analyticity.)
Let $R_\rho$ denote the $R$-transform
of the Bernoulli law $\lam_\rho:=(\delta_{-\rho}+\delta_{+\rho})/2$, that is,
$$R_\rho(z)=\frac{\sqrt{1+4\rho^2 z^2}-1}{2\rho 
z}=\frac{2z\rho}{\sqrt{1+4\rho^2z^2}+1}\,,$$
see \cite[Definition 5.3.22 and Exercise 5.3.27]{AGZ},
so that we have
\begin{equation}\label{eq333}
G_U(z_1,z_2)
=\frac{1}{2}G(z_1,z_2)R_\rho(G(z_1))\,.
\end{equation}
Repeating the computation 
with $G_{U^*}$, we have
$G_{U^*}=G_U$. Algebraic manipulations yield
\begin{eqnarray}
G_T(z_1,z_2)&=&z_2G(z_1,z_2)-G(z_1)\,,\label{eq1}\\
2\rho G_U(z_1,z_2)+G_T(z_1,z_2)&=& z_1 G(z_1,z_2)-G_T(z_2)\label{eq2}\,.
\end{eqnarray}
Therefore, we get by substituting 
\eqref{eq333} and \eqref{eq1} into \eqref{eq2}
that 
\begin{equation}\label{eq3a}
\rho G(z_1,z_2)R_\rho(G(z_1)) +z_2 G(z_1,z_2)-G(z_1)
=z_1 G(z_1,z_2)-G_T(z_2)\,,
\end{equation}
which in turns gives, for any $z_1,z_2\in {\mathbb C}^+$, 
\begin{equation}\label{eq4}
G(z_1,z_2)\left( \rho R_\rho(G(z_1))+z_2-z_1
\right)=G(z_1)-G_T(z_2)\,.
\end{equation}
Thus, 
\begin{equation}\label{eq6}
	G_T(z_2)=G(z_1) \quad \mbox{\rm when \ } z_2=z_1 -\rho R_\rho
	(G(z_1))\,.
\end{equation} 
The choice of $z_2$ as in \eqref{eq6} is allowed for any $z_1\in
\C^+$ because $G:\C^+\to \C^-$ and we can see that 
$R:\C^-\to\C^-$. Thus $\Im (z_2)\ge \Im (z_1)>0$, implying that
such $z_2$ belongs to the domain of
$G_T$.

The relation \eqref{eq6} is the Schwinger--Dyson equation in our
setup. It gives an implicit equation for $G(\cdot)$ in terms of $G_T(\cdot)$.
Further, for $z$ with large modulus,
$G(z)$ is small and thus $z\mapsto z-\rho R_\rho
(G(z))$ possesses a non-vanishing
derivative, and further is close to $z$. Because
$G_T$ is analytic in the upper half plane and its derivative
behaves like $1/z^2$
at infinity, it follows by the implicit function 
theorem that \eqref{eq6} uniquely determines 
$G(\cdot)$ in a neighborhood of $\infty$. By analyticity, it thus fixes 
$G(\cdot)$ in the upper half plane (and in fact, everywhere except in
a compact subset of $\mathbb{R}$), and thus determines uniquely the law
of ${\bf Y}$.

\begin{remark}
	\label{remark-law}
	Let $\mu_T$ denote the {\em spectral measure} of
	$T$, that is $\int fd\mu_T=\mu(f(T))$ for any $f\in C_b(\R)$.
	We emphasize that $G_T$ is {\rm not} the Stieltjes transform
	of 
	$\mu_T$; rather, it is the Stieltjes transform
	of the symmetrized version of the law of $T$, that is
	of the probability measure $\tilde \mu_T$.
	 With this convention, \eqref{eq6} is equivalent to the
	 statement that 
the law of ${\bf Y}$, denoted $\mu_Y$, equals the
{\em free convolution} of $\tilde \mu_T$ and 
$\lam_\rho$, i.e. 
$\mu_Y= \tilde \mu_T\boxplus \lam_\rho$. 
\end{remark}
\begin{remark}
	\label{rem-stransf}
	We provide, following \cite{haageruplarsen},
an 
alternative characterization of $\mu_A$ and
its support. We first introduce some terminology from
\cite{haageruplarsen}.
		Consider a tracial non-commutative $W^*$-probability space
		$({\cal M},\tau)$. Let $u$ be Haar-distributed and
		let $h$ be a 
		$*$-free from $u$  self-adjoint element 
		(whose law will be taken to be $\Theta$).
		Let $\tilde\nu^z$ denote the  law of $|zI-uh|$.
		The {\it Brown measure} for $uh$ is defined as
		$$ \frac{1}{2\pi}\Delta_z\int \log |x|
		d\tilde\nu^z(x)\,,$$
		c.f. \cite[Pg. 333]{haageruplarsen}.
Recall that $\Theta(\{0\})=0$
by Assumption 
		\ref{eq-denbound}.
		By  \cite[Proposition 3.5]{haageruplarsen} and Remark
	\ref{remark-law} above, $\tilde \nu^z=\nu^z$, and therefore,
	$\mu_A$ in the statement of Theorem
	\ref{main-theo} is the Brown measure for $uh$.
	By \cite[Theorem 4.4 and Corollary 4.5]{haageruplarsen},
	the Brown measure $\mu_A$ is radially symmetric and possesses 
	a density $\rho_A$ that can be described as follows.
Let 
$\Theta^{\sharp 2}$ denote the push forward of $\Theta$ by the
map $z\mapsto z^2$, i.e. $\Theta^{\sharp 2}$ is the weak limit
of $\{L_{T_n^2}\}$. Let ${\cal S}$ denote the S-transform 
of $\Theta^{\sharp 2}$ (see \cite[Section 2]{haageruplarsen} for the
definition of the S-transform of a probability measure on $\R$
and its relation to the R-transform).  
Define 
$F(t)=1/\sqrt{{\cal S}(t-1)}$ on ${\cal D}=(0,1]$. Then, $F$
maps ${\cal D}$ to the interval 
$$(a,b]=\left(\frac{1}{(\int x^{-2} d\Theta(x))^{1/2}},
\left( \int x^2 d\Theta(x)\right)^{1/2}\right]\,,
$$ and 
has an analytic continuation to a neighborhood 
of ${\cal D}$, and $F'>0$ on ${\cal D}$.
Further, with $\mu_A$ as above, $\rho_A(re^{i\theta})=\rho_A(r)$
and it holds that
\begin{equation}
	\label{eq-radial}
	\rho_A(r)=\left\{\begin{array}{ll}
		\frac{1}{2\pi r F'(F^{-1}(r))}\,,& r\in (a,b]\,,\\
		0\,,& otherwise.
	\end{array}
	\right.
\end{equation}
Finally, $\rho_A$ has an analytic continuation to 
a neighborhood of $(a,b]$, and $\mu_A$ is a probability measure, 
see \cite[Pg 333]{haageruplarsen}.
\end{remark}

In the next section, we will need the following estimate.
\begin{lemma}
	\label{lem-apriori}
	If $|\Im G_T(\cdot)|\leq \kappa_1$ on $\{z: \Im(z)\geq \epsilon\}$ then 
	$|\Im G(\cdot)|\leq \kappa_1$ on $\{z: \Im(z)\geq \epsilon\}$.
\end{lemma}
{\bf Proof} Recall that if $z\in \C^+$ then $G(z)\in \C^-$ and
also $R_\rho(G(z))\in \C^{-}$ because $R_\rho$ 
maps $\C^-$ into $\C^-$ (regardless
of the branch of the square root taken at each point). 
Thus, $y=z-R(G(z))$ has $\Im(y)\geq \Im(z)$. Therefore,
if $\Im(z)\geq \epsilon$ then
$|\Im G(z)|=|\Im G_T(y)|\leq \kappa_1.$ \qed

\subsection{Finite $n$ equations and convergence}
\label{subsec-finiteneq}
We next turn to the evaluation of the law of ${\bf Y}_n$. We assume
throughout that the sequence $T_n$ is uniformly bounded by some constant
$M$, that $L_{T_n}\to \mu_T$ weakly in probability,
and further that 
\eqref{eq-denbound} 
is satisfied.
All constants in this section are independent of $\rho$, but depend 
implicitly on $M$, the uniform bound on $\|T_n\|$ and on $\rho$.

 Recall first,
see \cite[(5.4.29)]{AGZ}, that by invariance of the Haar measure under 
unitary
conjugation, with $P\in\mathbb{C}\langle T,U,U^*\rangle$
a noncommutative polynomial (or a product of Stieltjes functionals),
\begin{equation}\label{master1}
E[\frac{1}{2n}\tr\otimes\frac{1}{2n}\tr(\partial P ({\bf T}_n, {\bf U}_n,
{\bf U}_n^*))]=0\,.
\end{equation}
This key equality can be proved by noticing that 
 for any $n\times n$ matrix $B$ such that $B^*=-B$, for
any $(k,\ell)\in  [1,n]$, if we let $U_n(t)= U_n e^{tB}$ and
construct ${\bf U_n}(t)$ and ${\bf U_n}^*(t)$ with this unitary matrix,
\begin{equation}\label{master11}
0=\partial_t E[ \left( P({\bf T}_n, {\bf U}_n (t),
 {\bf U}_n^*(t))\right)_{k,\ell}]
= E[ \left(\partial  P ({\bf T}_n, {\bf U}_n,
{\bf U}_n^*)\sharp {\bf B}\right)_{k,\ell}] 
\end{equation}
with ${\bf B}=\left(\begin{array}{ll}
0&0\\
0&B\\
\end{array}\right)$.
Letting $\Delta({k,\ell})$ be the $n\times n$ matrix so that
$\Delta({k,\ell})_{i,j}=1_{i=k}1_{j=\ell}$, we can choose in the last equality
$B= \Delta({k,\ell})-\Delta({\ell,k})$ or\\
$B=i\left(\Delta({k,\ell})+\Delta({\ell,k})\right)$. Summing the two
resulting
equalities and then summing over $k$ and $\ell$ yields
\eqref{master1}.

We denote by $G^n$ the quantities as defined in 
\eqref{eq-070809a},
but with $E[\frac{1}{2n}\tr]$ replacing $\mu$ and the superscript or subscript
$n$ attached
to all variables, so that for instance
$$G^n(z)=E[\frac{1}{2n}\tr\left((z-{\bf Y}_n)^{-1}\right)]\, .$$ 
We get  by taking $P=(z_1-{\bf Y}_n)^{-1}(z_2-{\bf T}_n)^{-1} {\bf U}_n$
that
\begin{equation}\label{eqUN}
\frac{1}{2}G_U^n(z_1,z_2)=-\rho
G_{U}^n(z_1,z_2)G_U^n(z_1)+\frac{\rho
}{4} G^n(z_1,z_2)
G^n(z_1)+O(n,z_1,z_2)\,,
\end{equation}
with 
$$O(n,z_1,z_2)=E\left[(\frac{1}{2n}\tr- E[\frac{1}{2n}\tr])\otimes 
(\frac{1}{2n}\tr- E[\frac{1}{2n}\tr])\partial
(z_1-{\bf Y}_n)^{-1}(z_2- {\bf T}_n)^{-1}{\bf U}_n \right]\,.$$
Further, by the standard concentration inequality for ${\cal H}_n$, see
\cite[Corollary 4.4.30]{AGZ}, for any smooth function $P:{\mathcal U}(n)\rightarrow \C$,
\begin{equation}
	\label{eq-conc}
	\given E \left[
\left(\frac{1}{2n}\tr (P)- 
E[\frac{1}{2n}\tr](P)\right)^2\right] \given
\le \frac{1}{n^ 2}\|P\|_L^2\,,
\end{equation}  
with $\|P\|_L$ the Lipschitz constant of $ P$ given by
$$\|P\|_L= \| D P\|_\infty$$ if
$D$ is the cyclic derivative given by $D=m\circ \partial$ 
with $m(A\otimes B)=BA$ and 
$\|DP\|_\infty$ denotes the operator norm. 
(The appearance of the cyclic derivative in the evaluation of the Lipshitz constant can be seen by approximating $P$ by polynomials.)
Applying \eqref{eq-conc} to each 
term of  $\partial P$  (recall  formula \eqref{caf}),
we get that for $\Im(z_1),\Im(z_2)>0$, 
and with $a\wedge b=\min(a,b)$,
$$| O(n,z_1,z_2)|\le 
\frac{C\rho^2}{n^2 |\Im(z_2)| \Im(z_1)^2 (\Im(z_1)\wedge 1)}\,.$$
(The inequality uses that for any Hermitian matrix, $\|(z-H)^{-1}\|_\infty
\leq 1/|\Im(z)|$.)
Multiplying by $z_2$ and taking the limit as $z_2\to\infty$
we deduce from \eqref{eqUN} that 
\begin{equation}
	\label{eq-120709a}
	\rho (G^n(z_1))^2=2G^n_U(z_1)(1+2\rho G^n_U(z_1))-O_1(n,z_1)\,,
\end{equation}
where 
\begin{eqnarray*} 
O_1(n,z_1)&=&4 E\left[(\frac{1}{2n}\tr- E[\frac{1}{2n}\tr])\otimes 
(\frac{1}{2n}\tr- E[\frac{1}{2n}\tr])\partial
(z_1-{\bf Y}_n)^{-1}{\bf U}_n \right]\\
&=&O\left(\frac{\rho^2}{n^2\Im(z_1)^2 (\Im(z_1)\wedge 1)}\right)\,.
\end{eqnarray*}
In particular, 
\begin{equation}
	\label{eq-070809d}
	G^n_U(z_1)=\frac{1}{4\rho}(-1+\sqrt{1+4\rho^2
	G^n(z_1)^2+4O_1(n,z_1)})\,,
\end{equation}
with again the choice of the square root determined by analyticity 
and behavior at infinity.

Recalling that
\eqref{eq1} and \eqref{eq2} remain true when we add the subscript $n$
and  combining these with \eqref{eqUN},
we get
\begin{equation}\label{mastereq}
G^n(z_1,z_2)\left( \frac{ \rho^2 G^n(z_1)}{(1+2\rho G_U^n(z_1))}+z_2
-z_1\right)= G^n(z_1)- G_{T_n}(z_2)
+\tilde O(n,z_1,z_2)\,,\end{equation}
with 
$$\tilde O(n,z_1,z_2)=\frac{2 O(n,z_1,z_2)}{(1+2\rho G_U^n(z_1))}\,.$$
Hence, 
if  we define
\begin{equation}
	\label{eq-080709c}
	z_2=\psi_n(z_1):=z_1-\frac{\rho^2 G^n(z_1)}{(1+2\rho 
	G^n_U(z_1))}\,,
\end{equation}
then
$$G^n(z_1)=G_{T_n}(z_2)-\tilde O(n,z_1,z_2)\,,$$
and therefore
\begin{equation}\label{cvp0}
G^n(z_1)=
G_{T_n}(\psi_n(z_1))-\tilde O(n,z_1,\psi_n(z_1))\,.
\end{equation}
Equation \eqref{cvp0} holds at least when $\Im( z_2)>0$ for
$z_2$ as in \eqref{eq-080709c}.
In particular, for 
$\Im(z_1)$ 
large (say larger than some
$M$), it holds that $G^n(z_1)$ and $G_U^n(z_1)$ are small,
implying that $z_2$ is well defined with $\Im(z_2)>0$.
Assume $L_{T_n}$ converges towards $L_T$ so that $G_{T_n}$ converges
to $G_T$ on $\mathbb C^+$. Then, 
 the limit points of the sequence of  uniformly continuous  functions
$(G^n(z),G_U^n(z))$ on $\{z:\Im (z)\ge M\}$ 
satisfy \eqref{eq3} and \eqref{eq6} and therefore equal
$(G(z),G_U(z))$  on $\{z:\Im (z)\ge M\}$ 
by uniqueness of the solutions to these equations. Hence,
taking $n\to\infty$
then implies that $G^n\to G$ in a neighborhood
in the upper half plane close to $\infty$. Since $G^n$ and $G$ are
Stieltjes transforms of probability measures,  
we have now shown the following (see Remark~\ref{remark-law}). 
\begin{lemma}
	\label{lem-S-conv}
	Assume $L_{T_n}$ converges weakly in probability to a compactly
	supported probability measure $\mu_T$. Then,
	$L_{{\bf Y}_n}$ converges weakly, in probability, to
	$\mu_Y= \tilde \mu_T\boxplus \lam_\rho$. In particular, if 
	$L_{T_n}$ converges weakly in probability to
	a probability measure $\Theta$, then
	for any
	$z\in \C$,  $\nu_n^z$ converges weakly in probability
	to $\tilde \Theta \boxplus \lam_{|z|}$.
\end{lemma}
\noindent
(Recall that
$\tilde \Theta$
is the symmetrized version of $\Theta$ and
note that for $z=0$, the statement of the lemma is trivial.)

Lemma \ref{lem-S-conv} completes the proof of Step one in our program. 
To be able to complete Step two, we need to obtain quantitative information
from the (finite $n$) Schwinger--Dyson equations \eqref{cvp0}: 
our goal is to show that the left side remains bounded in a domain of 
the form $\{z\in \mathbb C^+: \Im (z)>n^{-c}\}$ for some $c>0$.
Toward this end, we will show that in such a region,
$\psi_n$ is analytic, 
$\Im \psi_n(z)>(\Im (z)/2)\wedge C$ for some positive constant $C$
and $\tilde O(n,z_1,\psi_n(z_1))$ is analytic and bounded there.
This will imply that \eqref{cvp0} extends by analyticity to this region,
and our assumption on the boundedness of $G_{T_n}$ will lead to the conclusion.

As a preliminary step, note that
$G^n(\cdot)$ and $G^n_U(\cdot)$ are analytic in $\mathbb C^+$.
We have the following.
\begin{lemma}
	\label{lem-130709a}
	There exist  constants $C_1,C_2$ such that for all
	$z\in \C^+$ with $\Im(z)>C_1 n^{-1/3}$ and all $n$
	large, it holds that
	\begin{equation}
		\label{eq-071209c}
		|1+2\rho G_U^n(z)|>C_2\rho [\Im(z)^3 \wedge 1]\,.
	\end{equation}
\end{lemma}
\noindent
	{\bf Proof}
	Since $G_U^n(z)$ is asymptotic to 
	$1/z^2$ at infinity, we may and
	will restrict attention to some fixed ball $B_R\subset \C$, whose interior contains the support of ${\bf Y}$. But
	$$\Im(G^n(z))=-\Im(z)\int \frac{d\mu_{{\bf Y}_n}(x)}{
	(\Re(z)-x)^2+\Im(z)^2}$$
	and therefore, as $(\Re(z)-x)^2+\Im(z)^2\le 4 R^2$ for
all $z,x\in B(0,R)$
	\begin{equation}\label{270709}
|G^n(z)|\ge |\Im(G^n(z))|\geq \frac{|\Im(z)|}{4 R^2}\,.
\end{equation}
Moreover, 
since $|G_U^n(z)|\leq 1/|\Im(z)|$, 
	we deduce from 
	\eqref{eq-120709a} that
for some constant $c$ 
	independent of $n$ and all $n$ large,  
	$$|G^n(z)|^2\le  \frac{2|1+2\rho G_U^n(z)|}{\rho
	|\Im(z)|}+\frac{c\rho}{n^2 \Im(z)^2
(\Im(z)\wedge 1) }\,.$$
Combining this estimate and \eqref{270709},  we get that
\begin{equation}
	\label{eq-071309d}
	\frac{2|1+2\rho G_U^n(z)|}{\rho
	|\Im(z)|}\geq \frac{ |\Im(z)|^2}{16 R^4}-
\frac{c\rho }{n^2 \Im(z)^2 (\Im(z)\wedge 1) }\geq 
	 \frac{|\Im(z)|^2}{32R^4}
\,,
\end{equation}
as soon as $\Im(z)>C_1n^{-1/3}$ for an appropriate $C_1$, and
$|z|<R$. The conclusion follows.

\qed

As a consequence of Lemma 
	\ref{lem-130709a} and the analyticity of $G^n$ and $G_U^n$
	in $\C^+$, we conclude that $\psi_n$ is
	analytic in $\{z: \Im(z)>C_1n^{-1/3}\}$, for all $n$ large.

Our next goal is  to check the analyticity of
$z\ra \tilde O(n,z,\psi_n(z))$ for $z\in \C^+$ with imaginary part
bounded away from $0$ by a polynomially decaying (in $n$) factor.
Toward
this end, 
we now verify that
$\psi_n(z)\in\mathbb C^+$ for $z$ up to a small distance from
the real axis. 
\begin{lemma}
	\label{lem-130709b}
	There exists a constant $C_3$ such that if $\Im(z)>C_3
n^{-1/4}$,
	then
$\Im(\psi_n(z))\ge \Im(z)/2.$
\end{lemma}
\noindent
{\bf Proof}
Again, 
because both $G^n(z)$ and $G^n_U(z_1)$ tend to $0$ at infinity,
we may and will restrict attention
to $\Im(z)\leq R$ for some fixed $R$. 
We divide the proof to two cases, as follows.
Let $\e_n=n^{-1/2}$, and set
$\Delta_n=\{z\in \C^+: 
|\rho G^n(z) +i/2|\ge \e_n\}$.

Then, for any $z\in \Delta_n$, and whatever choice of branch of
the square root made in \eqref{eq-070809d},
if $\e_n^{-1/2}O_1(n,z)$ is small enough (smaller than $\e_n/2$ is fine), then
that choice can be extended to include a neighborhood of
the point $w=G^n(z)$ such that
with this choice, the function
$r_\rho (w)=\frac{1}{4\rho}(-1+\sqrt{1+4\rho^2w^2})$ 
is Lipschitz in the sense that
\begin{equation}
	\label{eq-120709f}
	|G^n_U(z)-r(G^n(z))|\le C \e_n^{-\frac{1}{2}}  O_1(n,z)/\rho\,.
\end{equation}
On the other hand, again from 
\eqref{eq-120709a},
$$ \left|\frac{\rho G^n(z)}{1+2\rho G^n_U(z)}-
\frac{2G_U^n(z)}{G^n(z)}\right|\leq
C\frac{|O_1(n,z)|}{|G^n(z)(1+2\rho G_U^n(z))|}\,.$$
Combining the last display with the relation $R_\rho(\theta)=
2r_\rho(\theta)/\theta$,
 \eqref{eq-120709f} and \eqref{270709},
one obtains that for $z\in \Delta_n$,
\begin{eqnarray*}
	\left|\frac{\rho G^n(z)}{1+2\rho G^n_U(z)}- \rho R_\rho
	(G^n(z))
	\right|&\leq&\left|\frac{2r(G^n(z))}{G^n(z)}-
\frac{2G_U^n(z)}{G^n(z)}\right|+\left|\frac{\rho G^n(z)}{1+2\rho G^n_U(z)}-
\frac{2G_U^n(z)}{G^n(z)}\right|\nonumber\\
&\leq&
 C  \frac{| O_1(n,z)|}{\rho \e_n^{\frac{1}{2}} |G^n(z)|}
+C\frac{|O_1(n,z)|}{|G^n(z)(1+2\rho G^n_U(z))|}
\nonumber\\
&\le& C\frac{|O_1(n,z)|}{\rho \e_n^{1/2}|\Im(z)|}
	+C\frac{|O_1(n,z)|}{\rho \Im(z)^4}\nonumber\\
	&\leq& \frac{C\rho }{n^2 |\Im(z)|^4}\left(\frac{1}{\e_n^{1/2}}+
	\frac{1}{|\Im(z)|^3}\right)\nonumber\\
	&\leq&
	 \frac{C\rho }{n^2 |\Im(z)|^4}\left({n^{1/4}}+
	\frac{1}{|\Im(z)|^3}\right)
	\,.
\end{eqnarray*}
Since the above  right hand side is smaller than $\Im(z)/2$ 
for $\Im(z)>n^{-1/4}$,
 we conclude
 that for $z\in\Delta_n\cap\{\Im(z)>n^{-1/4}\}$
 \begin{equation}
	 \label{eq-071309h}
\Im \left(\frac{\rho G^n(z)}{1+2\rho G^n_U(z)}\right)\le \frac{1}{2}\Im (z)
\end{equation}
as, regardless of the branch taken in the definition of $R_\rho(\cdot)$,
$\Im
R_\rho(G^n(z))\le 0$.

On the other hand, when
$z\in \C^+\setminus\Delta_n$ and $\Im(z)>n^{-1/4}$,
then we have from
\eqref{eq-070809d} that for all $n$ large,
$$ |\rho G_U^n(z)+1/4|\leq \frac12\sqrt{\e_n+|O_1(n,z)|}\leq \frac18\,.$$
Thus, under these conditions,
\begin{eqnarray*}
\Im\left(\frac{\rho G^n(z)}{1+2\rho G^n_U(z)}\right)
&=& 
\Im\left(\frac{2\rho  G^n(z)}{1+4(\rho G^n_U(z)+1/4)}\right)
\\&\leq& 
2\rho \Im( G^n(z))
+16\rho |\Im( G^n(z))|
|\rho G_U^n(z)+1/4|\,,
\end{eqnarray*}
where we used that for $|a|\leq 1/2$, $|a/(1-a)|\leq 2|a|$.
Consequently, since 
$\rho G^n(z)$ is uniformly
bounded on 
$\C^+\setminus \Delta_n$ and $\Im(G_n(z))<0$ there, 
we get
$$\Im\left(\frac{\rho G^n(z)}{1+2\rho G^n_U(z)}\right)
\leq 
C\sqrt{\e_n+|O_1(n,z)|}\leq 
 Cn^{-1/4}\,.$$
We thus conclude from the last display
and \eqref{eq-071309h} the existence of 
a constant $C_3$ such that if $\Im(z)>C_3n^{-1/4}$ then
$$\Im(\psi_n(z))=\Im(z)-
\Im\left(\frac{\rho G^n(z)}{1+2\rho G^n_U(z)}\right)\geq \Im(z)/2\,,$$
as claimed.
\qed

From Lemma \ref{lem-130709b} we 
thus conclude the analyticity of
$z\ra \tilde O(n,z,\psi_n(z))$ in $\{z: \Im(z)\geq C_3 n^{-1/4}\}$, and
thus, due to 
\eqref{eq-080709c} and \eqref{cvp0}, 
$ \rho G^n(z)/(1+2\rho G_U^n(z))$  is
also analytic there (compare with Lemma \ref{lem-130709a}). 
In particular, the equality 
\eqref{cvp0}
extends by analyticity to this region.

We have made all preparatory steps
in order to state the main result of this subsection.
\begin{lemma}
	\label{lem-tildegest}
	There exist positive finite
  constants $C_6,C_7,C_8$ such that, for  
	$n>C_6$ and all $z\in {\cal E}_n:=\{z: \Im(z)>n^{-C_{7}}\}$,
	\begin{equation}
		\label{eq-150709b}
		|\Im G^n(z)|\leq C_8\,.
	\end{equation}
\end{lemma}

\noindent
{\bf Proof}
This is immediate from Lemma 
 \ref{lem-130709a}, Lemma \ref{lem-130709b},
 the definition of $\psi_n$, the assumption 
 \eqref{eq-denbound}
 on $G_{T_n}$,
 and the equality \eqref{cvp0}.
 \qed
\section{Tail estimates for $\nu_n^z$}
\label{sec-tailnu}
For $R>0$, let 
$ B_R=\{z\in \C: |z|\in [0,R]\}$.
Our goal in this short section is to 
prove the following proposition.
\begin{proposition}
	\label{prop-lowertail}
	(i) 
	Under the assumptions of Theorem \ref{main-theo},
for Lebesgue almost every $z\in \C$,
\begin{equation}
	\label{eq-190709a}\lim_{\epsilon\downarrow 0}
	\limsup_{n\to\infty}
	E[1_{{\mathcal G}_n}\int_{0}^{\epsilon}\log|x| d\nu_n^z(x)]=0
\,.
\end{equation}
Consequently, for Lebesgue any  $z\in \mathbb C$,
\begin{equation}
	\label{eq-190709b}
	\int\log|x| d\nu_n^z(x)\to
	\int\log|x| d\nu^z(x)\,,
\end{equation}
	in probability.\\
	(ii) Fix $R>0$. For any smooth compactly supported deterministic
	function $\phi$ on 
	$ B_R$,
\begin{equation}
	\label{eq-280709a}
	\int \phi(z) \int\log|x| d\nu_n^z(x) dm(z)\to
	\int\phi(z) \int\log|x| d\nu^z(x)dm(z)\,,
\end{equation}
in probability.
\end{proposition}
Before bringing the proof of Proposition \ref{prop-lowertail}, 
we recall the following elementary lemma.
\begin{lemma}
	\label{lem-stap}
	Let $\mu$ be a probability measure on $\R$. For any real
	$y>0$, 
	it holds that
	\begin{equation}
		\label{eq-mubound}
		\mu((-y,y))\leq 2y |\Im G(iy)|\,.
	\end{equation}
\end{lemma}
\noindent
{\bf Proof} We have
$$-\Im(G(iy))=\int\frac{y}{y^2+x^2} \mu(dx)\geq \int_{-y}^y 
\frac{y}{y^2+x^2} \mu(dx)\geq \frac{1}{2y} \mu((-y,y))\,,$$
from which \eqref{eq-mubound} follows. \qed

We can now provide the\\
{\bf Proof of Proposition 
\ref{prop-lowertail}} 
	 
(i) Assume $z\in  B_R$ for some $R>0$.
  By 
\eqref{eq-300609b}, we can replace the lower limit of integration
in 
\eqref{eq-190709a}
with $n^{-\delta}$. Let $G_n^z$ denote the Stieltjes transform of
$E[\nu_n^z]$. 
By Lemma \ref{lem-tildegest} and 
Lemma \ref{lem-apriori}, there exist positive
constants $c_1=c_1(R),c_2=c_2(R)$ such that
whenever $\Im(u)>n^{-c_1}$, it holds that
$|\Im G_n^z(u)|<c_2$.  We may and will assume that $c_1<\delta$.

Since $G_n^z$ is the Stieltjes transform of $E[\nu_n^z]$,
 by Lemma \ref{lem-stap}, we have for any $y>0$ that
$$E[\nu_n^z((-y,y))]\le E[\nu_n^z((-y\vee n^{-c_1},y\vee n^{-c_1}))]
\le 2 c_2y\vee n^{-c_1}\, .$$

 Thus, we get that for any $z\in  B_R$ and with $\alpha\in [1,2]$, 
\begin{eqnarray*}
        &&
	E[\int_{n^{-\delta}}^{\epsilon}(|\log x|)^\alpha d\nu_n^z(x)]
        \\
        &\leq &E[
	 \int_{n^{-\delta}}^{n^{-c_1}} (|\log x |)^{\alpha}
d\nu_n^z(x) + \int_{n^{-c_1}}^{\epsilon}(|\log  x |)^\alpha
d\nu_n^z(x)]\\
&\leq &
        (  \delta\log n)^\alpha E[\nu_n^z((-n^{-c_1},n^{-c_1}))]
        \\
	&&
	+\sum_{j=0}^{J}E[ \nu_n^z((-2^{(j+1)} n^{-c_1}, 2^{(j+1)}
n^{-c_1}))](\log
        (2^j n^{-c_1}))^\alpha\,,
\end{eqnarray*}
	where $2^{J-1} n^{-c_1}<\epsilon\leq 2^Jn^{-c_1}$. Note that
	by Lemma 
\ref{lem-stap} and the estimate on $G_n^z$, for $j\geq 0$,
$$E[\nu_n^z( (-2^jn^{-c_1},2^jn^{-c_1}))]\leq 2^{j+1}c_2 n^{-c_1}\,.$$
We conclude that
\begin{equation}
	\label{eq-280709b}
	E[\int_{n^{-\delta}}^{\epsilon}|\log x|^{\alpha} d\nu_n^z(x)]
	\leq C \epsilon |\log(\epsilon)|^\alpha\,,
\end{equation}
	where the constant $C=C(R)$.
	To obtain the estimate \eqref{eq-190709a}, we will consider
	 $\alpha=1$ and argue as follows.
Due to \eqref{eq-300609b}, for $\alpha< 2$ 
we have
\begin{eqnarray*} 
& &	E[{\bf 1}_{{\cal G}_n}
\int_{0}^{n^{-\delta}}|\log x|^{\alpha} d\nu_n^z(x)] \\
&\le& 	E[{\bf 1}_{{\cal G}_n}\nu_n^z([-n^{-\delta}, n^{-\delta}])
		{\bf 1}_{\{\sigma_n^z<n^{-\delta}\}}
		|\log \sigma_n^z|^\alpha]\\
 &\le& E\left[\left(\nu_n^z([-n^{-\delta},
n^{-\delta}])\right)^{\frac{2}{2-\alpha}}\right]^{\frac{2-\alpha}{2}}
 E[{\bf 1}_{{\cal G}_n}
		{\bf 1}_{\{\sigma_n^z<n^{-\delta}\}}
		|\log \sigma_n^z|^2]^{\frac{\alpha}{2}}\\
\end{eqnarray*}
by H\"{o}lder's inequality. The first factor goes to zero because
$$ E\left[\left(\nu_n^z([-n^{-\delta},
n^{-\delta}])\right)^{\frac{2}{2-\alpha}}\right] \le E\left[\nu_n^z([-n^{-\delta},
n^{-\delta}])\right]\le
2c_{2} n^{-c_1}.$$ 
By \eqref{eq-300609b},
the second factor is bounded by $(\delta')^{\alpha/2}$. 
 We thus get \eqref{eq-190709a}
from \eqref{eq-280709b}.  
By Chebycheff's inequality, the convergence
in expectation implies the convergence in probability and therefore
for any $\delta,\delta'>0$ there exists $\epsilon>0$ small enough so that
$$\lim_{n\ra\infty} P(\int_0^\epsilon |\log x|  d\nu_n^z(x)>\delta)<\delta'$$
On the other hand, $\int_\epsilon^\infty \log|x|  d\nu_n^z(x)$
converges to 
 $\int_{\epsilon}^{\infty} \log|x|  d\nu^z(x)$ 
by the weak convergence of $\nu_n^z$ to $\nu^z$ in probability for any
$\epsilon>0$, and $\int_0^\epsilon \log |x| d\nu^z(x)$ converges to $0$ as
$\epsilon\to 0$ since $\nu^z$ has a bounded density 
by Lemma \ref{lem-apriori}. Hence, we get 
	 \eqref{eq-190709b}.

	 (ii) 
	 Define the functions 
	 $f_n^i: B_R\to \R$, $i=1,2$  by	 \begin{eqnarray*}
		 f_n^1(z)
		 &=&{\bf 1}_{\GG_n}
		 {\bf 1}_{\|T_n\|\leq M} 
		 \int_0^{n^{-\delta}} \log (x)d\nu_n^z(x)\,,\\
		 f_n^2(z)
&=&{\bf 1}_{\GG_n}
		 {\bf 1}_{\|T_n\|\leq M}
\int_{n^{-\delta}}^\infty \log (x)d\nu_n^z(x)\,,
\end{eqnarray*}
and set $f_n(z)=f_n^1(z)+f_n^2(z)$. Because 
$\nu_n^z$ is
supported in $B_{R+M}$ on $\|T_n\|\le M$ for all $z\in  B_R$, 
$f_n$ is bounded above
by $\log(R+M)$. By
\eqref{eq-280709b}, $E[|f_n^2(\cdot)|^2$ is bounded, uniformly
in $z\in  B_{R}$. On the other hand, by 
\eqref{eq-300609b}, 
again uniformly in $z\in  B_R$,
$E(f_n^1(z)^2)<\delta'$, and therefore
$$E\int_{\tilde B_R} (f_n^1(z))^2 dm(z)<\infty\,.$$
Thus, $E\int_{\tilde B_R} |f_n(z)|^2 dm(z)<\infty$, and
in particular, the sequence of random variables 
$$\int_{\tilde B_R} \Big|{\bf 1}_{\GG_n}{\bf 1}_{\|T_n\|\leq M}
\int \log x d\nu_n^z(x)\Big|^2 dm(z)$$ is bounded in probability.
This uniform integrability
and the weak convergence 
	\eqref{eq-190709b} are enough to conclude, using
	dominated convergence (see
	\cite[Lemma 3.1]{taovuk} for a similar argument).
\qed

\section{Proof of Theorem \ref{main-theo}}
\label{sec-theoproof}
It clearly suffices to prove the theorem for deterministic
diagonal matrices $T_n$. (If $T_n$ is random, use the
independence of $(U_n,V_n)$ from $T_n$ to apply the deterministic version,
after restricting attention to matrices $T_n$ belonging to a set whose
probability approaches $1$).
By Proposition \ref{prop-lowertail}, see 
        \eqref{eq-280709a},
we have, with
$h(z):= \int \log|x| d\nu^{z}(x)$,
         that
for any $R$ and any smooth function $\psi$ on $\tilde B_R$,
 $$\int \psi(z) dL_{A_{n}}(z)  \rightarrow 
\frac{1}{2\pi} 
 \intt_{\C} \Delta \psi (z) \; h(z) dm(z)\,,$$ 
in probability. Since the sequence 
$L_{A_n}$ is tight,
it thus follows that it converges, in the sense of distribution,
to  the  measure
$$\mu_A:=\frac{1}{2\pi} 
 \Delta_z  h(z) \,.$$ 
From 
Remark 
\ref{rem-stransf} (based on
\cite[Corollary 4.5]{haageruplarsen}), 
we have that $\mu_A$ is a probability measure that
possesses a radially symmetric density
$\rho_A$ satisfying the properties stated in parts b and c of the
theorem.
\qed

\section{Proof of Theorem \ref{cor-FZ}}
\label{sec-corproof}
We let $X_n$ be as in the statement of the corollary and write
$X_n=P_n T_n Q_n$ with $P_n,Q_n$ unitary and $T_n$ diagonal with entries
equal to the 
singular values $\{\sigma_i^n\}$ of $X_n$. Obviously,
 $\{P_n,Q_n\}_{n\geq 1}$ is a sequence of independent,
${\mathcal H}_n$-distributed matrices.
The joint distribution of the  entries of $T_n$
possesses a density on $\R_+^n$ which
is given
by the expression
$$ \tilde Z_n \prod_{i<j} |\sigma_i^2-\sigma_j^2|^2 
e^{-n \sum_{i=1}^n V(\sigma_i^2)}\prod_i \sigma_i d\sigma_i\,,$$ 
where $\tilde Z_n$ is a normalization factor,
see e.g. \cite[Proposition 4.1.3]{AGZ}. Therefore,
the squares of the singular values possess the joint density
$$ \hat Z_n \prod_{i<j} |x_i-x_j|^2 
e^{-n \sum_{i=1}^n V(x_i)}\prod_i  dx_i\,,$$ 
on $\R_+^n$.
In particular, it falls within the framework treated in
\cite{PS}. By  part (i) of Theorem 2.1 there,
there exist positive constants $M, C_{11}$ such that 
$P(\sigma_1> M-1)\leq e^{-C_{11} n}\,,$
and thus point 1 of the assumptions of Theorem \ref{main-theo}
holds. 
By equations \cite[ (2.26) and (2.27)]{PS}
and Chebycheff's inequality, we get that for $z$
with $\Im(z)>n^{-\kappa'}$
where   $\kappa<(1-\kappa')/2$, 
$$P\left( |G_{T_n}(z)-G_{\tilde\Theta}(z)|\ge \frac{1}{2\Im (z) n^{\kappa}}
\right)\le C |\Im (z)|^{-1} n^{2\kappa-1} \log n\,.$$
As the derivative of $ G_{T_n}-G_{\tilde\Theta}$ is bounded by
a constant multiple of $1/|\Im(z)|^2$, a covering argument and summation
shows that for $\kappa'<1/2$,
$$P\left(\sup_{z:|z|\le M\atop\Im (z) \ge n^{-\kappa'}}
 |G_{T_n}(z)-G_{\tilde\Theta}(z)|\ge \frac{1}{\Im (z) n^{\kappa}}\right)
\le M n^{4\kappa+2\kappa'-1}  \log n,$$
which  goes to zero for $\kappa\in (0,(1-2\kappa')/4)$.
Together with \cite[Equation (2.32)]{PS},
this proves point 3 of the assumptions.
Thus, it remains only to check point 2  of the assumptions.
Toward this end, 
define
${\cal G}_n=\{\sigma_1^n<M+1\}$ and note that
we may and will restrict attention to
$|z|<M+2$ when checking 
\eqref{eq-300609b}. We begin with the following proposition, 
due to \cite{SST}.
\begin{proposition}
	\label{prop-SST}
	Let $\bar A$ be an arbitrary $n$-by-$n$ matrix, and
	let $A=\bar A+\sigma N$ where $N$ is a 
	matrix with independent (complex) Gaussian 
	entries of zero mean and 
	 unit variances.
	Let $\sigma_n(A)$ denote the minimal singular value
	of $A$. Then, there exists a constant $C_{12}$
	independent of $\bar A$, $\sigma$ or $n$ such
	that
	\begin{equation}
		\label{eq-300709aa}
		P(\sigma_n(A)<x)\leq 
		C_{12}{n}\left(\frac{x}{\sigma}\right)^2\,.
	\end{equation}
\end{proposition}
The proof of Proposition \ref{prop-SST} is identical to 
\cite[Theorem 3.3]{SST}, with the required adaptation in
moving from real to complex entries. (Specifically, in the right side of 
the display in \cite[Lemma A.2]{SST}, $\epsilon
\sqrt{2/\pi}/\sigma$ is replaced by its square.) We omit further details.

On the event ${\cal G}_n$, 
all entries of the matrix
$X_n$ are bounded by a constant multiple of
$\sqrt{n}$. Let
$N_n$ be a Gaussian
matrix as in Proposition \ref{prop-SST}.
With $\alpha>2$ a constant to be determined below,
set
	$${\cal G}_n'=\{ \mbox{\rm all  entries of
	$n^{-\alpha/2}N_n$ are bounded by $1$
	}\}\,.$$ 
	Note that because $\alpha\geq 2$,
	on ${\cal G}_n'$, we have that
	$\sigma_1(n^{-\alpha}N_n)\leq 1$.
Define
$\bar A_n=zI-X_n$,
$\tilde A_n=\bar A_n+n^{-\alpha} N_n{\bf 1}_{ {\cal G}_n'}$ and
$A_n=
\bar A_n+n^{-\alpha} N_n$.
Then, by \eqref{eq-300709aa}, with $\sigma_n(A_n)$
denoting the minimal singular value
of $A_n$, we have
	\begin{equation}
		\label{eq-300709ab}
		P(\sigma_n(A_n)<x;{\cal G}_n)\leq 
		C_{12}x^2 n^{1+2\alpha}\,.
	\end{equation}
	If the estimate \eqref{eq-300709ab}
	  concerned $\bar A_n$ instead of $A_n$,
	 it would have been
	straightforward to check that
	point 2 of the assumptions of Theorem \ref{main-theo}
	holds (with an appropriately chosen
	$\delta$, which would depend on $\alpha$). 
	Our goal is thus to replace,
	in \eqref{eq-300709ab}, $A_n$ by $\bar A_n$, at the
	expense of not too severe degradation in the right side.
	This will be achieved in two steps: first, we will replace
	$A_n$ by $\tilde A_n$, and then we will
	 construct
	on the same probability space the matrix
	$X_n$ and a matrix $Y_n$ so that $Y_n$ is distributed
	like $X_n+n^{-\alpha} N_n{\bf 1}_{ {\cal G}_n'}$ 
	but $P(Y_n\neq X_n)$ is small.

	Turning to the construction, observe first that
		from \eqref{eq-300709ab},
	\begin{equation}
		\label{eq-300709abc}
		P(\sigma_n(\tilde A_n)<x;{\cal G}_n)\leq 
		C_{12}x^2 n^{1+2\alpha}+P( ({\cal G}_n')^c)
		\leq
		C_{12}[x^2 n^{1+2\alpha}+n^2 e^{-n^{\alpha}/2}
		]\,.
	\end{equation}
	Let $X_n^{(\alpha)}=X_n+n^{-\alpha}N_n{\bf 1}_{ {\cal G}_n'}$.
	Let $\{\theta_i\}$ and $\{\mu_i\}$ denote the eigenvalues 
	of $W_n=X_n X_n^*$ and of 
	$W_n^{(\alpha)}=(X_n^{ (\alpha)})(X_n^{ (\alpha)})^*$,
	respectively, arranged in decreasing order.
	Note that 
	the density of $X_n$ is of the form 
	$$ Z_n^{-1} e^{-n \tr(V({\bf x}{\bf x}^*))} d{\bf x}\,,$$
	where the variable ${\bf x}=\{x_{i,j}\}_{1\leq i,j\leq n}$ 
	is matrix valued and
	$d{\bf x}=\prod_{1\leq i,j\leq n} dx_{i,j}$,
	while that of $X_n^{(\alpha)}$ is of the form
	$$ Z_n^{-1} E_N[e^{-n \tr(V(
	({\bf x}+{\bf 1}_{{\cal G}_n'}
	n^{-\alpha}N_n)({\bf x}+{\bf 1}_{{\cal G}_n'}
	n^{-\alpha} N_n)^*))}] d{\bf x}\,,$$
	where $E_N$ denotes expectation with respect to the law of $N_n$, and
	$Z_n$ is the same in both expressions.
	Note that  $\sigma_1(X_n^{(\alpha)})\in [\sigma_1(X_n)-1,
\sigma_1(X_n)+1]$.
	Because $V(\cdot)$ is locally Lipschitz, we have
that if either $\sigma_1(X_n)\leq M+1$ or 
$\sigma_1(X_n^{(\alpha)})\leq M+1$, then
	there exists a constant $C_{13}$
	independent of $\alpha$ so that
	\begin{eqnarray*}
		|\tr(V(W_n)-V(W_n^{(\alpha)}))|&\leq&
		\sum_{i=1}^n |V(\theta_i)-V(\mu_i)|
		\leq C_{13} \sum_{i=1}^n|\theta_i-\mu_i|\\
		&\leq &
		C_{13} n^{1/2}
\left(\sum_{i=1}^n|\theta_i-\mu_i|^2\right)^{\frac 12}\\
		&\leq&
		C_{13}n^{1/2}\left( \tr( (W_n-W_n^{(\alpha)}
		)^2)\right)^{\frac 1 2 }\,,
	\end{eqnarray*}
	where the Cauchy--Schwarz inequality was used in the third inequality
	and the Hoffman--Wielandt inequality in the
	next  (see e.g. \cite[Lemma 2.1.19]{AGZ}).
	On the event ${\cal G}_n$, all entries of 
	$W_n-W_n^\alpha$ are bounded by $n^{(3-\alpha)/2}$.
	Therefore,
	\begin{equation}
		\label{eq-300709c}
		|\tr(V(W_n)-V(W_n^{(\alpha)}))|\leq
  n^{(C_{14}-\alpha)/2}\,,
  \end{equation}
where
	the constant $C_{14}$ does not depend on $\alpha$.
	In particular, if $\alpha>(C_{14}+1)\vee 2$ we obtain
	that  
	on ${\cal G}_n$, the ratio
	of the functions
	$f_n=e^{-n\tr(V( W_n))} $ and $g_n=e^{-n\tr(V(
W_n^{(\alpha)}))}$ is bounded e.g.
	by $1+n^{(C_{14}+1-\alpha)/2}$; in particular,
	it holds that
	\begin{eqnarray*}
		P(\sigma_1(X_n^{(\alpha)})< M)&\leq& (1+n^{(C_{14}+1-\alpha)/2})
	P(\sigma_1(X_n)< M)\\
	&\leq& (1+n^{(C_{14}+1-\alpha)/2})^2
	P(\sigma_1(X_n^{(\alpha)})< M)\,.
\end{eqnarray*}
	Therefore, the variational distance between 
	the law of $X_n$ conditioned on $\sigma_1(X_n)<M$
	and that of $X_n^{ (\alpha)}$ conditioned on
	$\sigma_1(X_n^{ (\alpha)})<M$,
	is bounded
	by
	$$ 
	4n^{(C_{14}+1-\alpha)/2}\,.$$

It follows that
one can construct a matrix $Y_n$ of 
law identical to the law of $X_n^{(\alpha)}$ conditioned on 
$\sigma_1(X_n^{\alpha})<M$,
together with $X_n$, on the same probability space so that
$$ P(X_n\neq Y_n; {\cal G}_n)\leq 4 n^{(C_{14}+1-\alpha)/2}
\leq n^{C_{15}-\alpha/2}\,.$$
Combined with
		\eqref{eq-300709abc}, we thus deduce that
		$$P(\sigma_n(\bar A_n)<x;{\cal G}_n)\leq 
		C_{12}x^2 n^{1+2\alpha}+n^{C_{16}-\alpha/2}
		\leq n^{C_{17}}x^{2/5}\,,
$$
where $\alpha$ was chosen as function of $x$.
This yields
immediately point 2 of the assumptions of Theorem \ref{main-theo}, if 
$\delta>5C_{17}/2$.

We have checked now that in the setup of Theorem 
\ref{cor-FZ}, all the assumptions of Theorem \ref{main-theo} hold.
Applying now the latter theorem completes the proof of Theorem
\ref{cor-FZ}.
\qed

\noindent
\begin{remark}
	The proof of Theorem \ref{cor-FZ} carries over to more
	general situations; indeed, $V$ does not need to be a polynomial,
	it is enough that its growth at infinity is
	polynomial and that it is locally Lipschitz, so that
	the results of \cite{PS} still apply. We omit further details.
\end{remark}

\section{Proof of Proposition \ref{cor-D} }
We take $T_n$ satisfying the assumptions
of Proposition \ref{cor-D} and consider $Y_n= U_n T_n  V_n +
n^{-\gamma} N_n$, with matrix of singular values $\tilde T_n$.
Note that $Y_n=\tilde U_n \tilde T_n \tilde V_n$
with $\tilde U_n,\tilde V_n$ following the Haar measure.
 We first show that $\tilde T_n$ 
also satisfies the assumptions of Theorem \ref{main-theo}
when  $\gamma>\frac 1 2$, except for the
second one. Since the singular values of $N_n$ 
follows the joint density
of Theorem \ref{cor-FZ} 
with $V(x)=\frac{1}{2}x^2$, it follows from the previous section that
$P(\|n^{-\frac{1}{2}} N_n\|>M)\le e^{-C_{11}n}$ and therefore
$\|\tilde T_n\|\le \|T_n\|+n^{-\gamma +\frac 1 2} \|n^{-\frac 1 2} N_n\|$ is bounded with overwhelming probability. 
Moreover, since $\tilde T_n=|T_n+n^{-\gamma} U_n^* N_n V_n^*|$,
on the event $\|N_n/\sqrt{n}\|\leq M$ we have
$$\left| G_{T_n}(z)-G_{\tilde T_n}(z)
\right|\le \frac{E[\|\tilde T_n-T_n\|{\bf 1}_{\|N_n/\sqrt{n}\|\leq 
M}]}{|\Im (z)|^2}\le 
\frac{C (\|T_n^{-1}\|,\|T_n\|)}{|\Im (z) |^2}n^{\frac{1}{2} -\gamma}
$$
with $C (\|T_n^{-1}\|,\|T_n\|)$ a finite constant depending only on
$\|T_n^{-1}\|,\|T_n\|$
which we assumed bounded. 
(In deriving the last estimate, we used that $\|(I+B)^{1/2}-I\|
\leq \|B\|$ when $\|B\|<1/2$.)
As a consequence, the third condition is satisfied since
$$\left| G_{\tilde\Theta}(z)-G_{\tilde T_n}(z)
\right|\le \frac{C (\|T_n^{-1}\|, \|T_n\|)}{|\Im (z) |^2}
n^{\frac{1}{2} -\gamma} +\frac{K}{n^\kappa |\Im (z)|}
\le \frac{K'}{n^{\gamma'} |\Im ( z)|}$$
with $\gamma'=\min\{\kappa, \frac{1}{2}(\gamma-\frac  12)\}$
and $\Im (z)\ge n^{-\max\{\frac{1}{2}(\gamma-\frac  12), \kappa'\}}$.
Hence, the results of Lemma \ref{lem-tildegest} hold and 
we need only 
check, as in Proposition \ref{prop-lowertail}, that
with $\nu_n^z$ the empirical measure of the 
singular values of $zI-Y_n$, 
	$$I_n:=E[1_{{\mathcal G}_n}\int_{0}^{n^{-\delta}}\log|x| d\nu_n^z(x)]$$
vanishes as $n$ goes to
infinity  for some $\delta>0$ and some set $\GG_n$ with overwhelming probability. But $\bar A_n=zI-Y_n= zI-  U_n T_n  V_n +n^{-\gamma}\tilde N_n$
with $\tilde N_n$ a Gaussian matrix, and therefore we can use
Proposition \ref{prop-SST} to obtain \eqref{eq-300709aa}
with $\sigma=n^{-\gamma}$, and the desired estimate on $I_n$.
\qed

{\bf Proof of Example \ref{exa}}
The first and the third hypotheses of Theorem \ref{main-theo}
are verified since $\mu$ is compactly 
supported and we assumed that the imaginary part of the Stieltjes transform
of its symmetrized version is uniformly bounded
on ${\mathbb C}^+$. For the third, note that if $F^{-1}$ is H\"older
continuous with
index $\alpha$, 
$$\left| G_{\tilde \Theta}(z)- G_{T_n}(z)\right|\le \sum_{i=1}^n
\frac{|s_{i+1}^n-s_i^n|
}{n|\Im (z)|^2}=\sum_{i=1}^n
\frac{| F^{-1} (\frac{i+1}{n})- F^{-1} (\frac{i}{n})|
}{n|\Im (z)|^2}\le C\frac{ n^{-\alpha}}{|\Im (z) |^2}$$
where we finally used that $F^{-1}$ is H\"older continuous with
index $\alpha$. \qed

\section{Extension to orthogonal conjugation}
In this section, we generalize 
Theorem \ref{main-theo} to the case where we conjugate 
$T_n$ by orthogonal
matrices instead of unitary matrices. 

\begin{theorem}\label{theo-orth}
Let $T_n$ be a sequence of diagonal matrices
satisfying the assumptions of Theorem \ref{main-theo}.
Let $O_n, \tilde O_n$ be two $n\times n$  independent 
matrices which  follow the Haar measure on the orthogonal group
and set $A_n=O_n T_n \tilde O_n$. Then, $L_{A_n}$ 
converges in probability to the probability measure $\mu_A$ described
in Theorem \ref{main-theo}.
\end{theorem}
{\bf Proof.} To prove the theorem, it is enough, following 
Section \ref{sec-theoproof}, to prove the analogue of Lemma \ref{lem-tildegest}
which in turn is based on the approximate Schwinger--Dyson equation
\eqref{mastereq} which is itself
a consequence of equation \eqref{master1}
and concentration inequalities.  To prove the analogue of
\eqref{master1} when $U_n$ follows the Haar measure 
on the orthogonal group,
observe that \eqref{master11} remains true with $B^t=-B$
which only leaves the choice $B=\Delta(k,\ell)-\Delta(\ell,k)$
possible. However,  taking this choice
and summing over 
 $k,\ell$, yields, if we denote $\tilde m(A\otimes B)=AB^t$,
$$E[\frac{1}{2n}\tr\otimes\frac{1}{2n}\tr(\partial P ({\bf T}_n, {\bf U}_n,
{\bf U}_n^*))]=\frac{1}{2n}
E[\frac{1}{2n}\tr\left((\tilde m\circ \partial P) ({\bf T}_n, {\bf U}_n,
{\bf U}_n^*)\right)].$$
The right hand side is small as $\tilde m\circ \partial P$ is uniformly
bounded. In fact, taking 
$P=(z_1-{\bf Y}_n)^{-1}(z_2-{\bf T}_n)^{-1} {\bf U}_n$, we find that
$\tilde 
m\circ \partial P$ is uniformly bounded by $2/(|\Im (z_2)|(|\Im (z_1)|\wedge 
1)^2)$  and therefore \eqref{eqUN} holds once we
add to  $O(n,z_1,z_2)$ the above right hand side which is at 
most of order $1/ n|\Im (z_2)|(|\Im (z_1)|\wedge 
1)^2$. Since our arguments did not require a very fine control
on the error term, we see that this change will not affect them.
Since concentration inequalities also hold under the
Haar measure on the orthogonal group, see \cite[Theorem 4.4.27]{AGZ}
and  \cite[Corollary 4.4.28]{AGZ}, the proof of Theorem
\ref{main-theo}
 can be adapted to
this set up. \qed

\section{Proof of Proposition \ref{corplus}}
We use again  Green's formula 
\begin{eqnarray*}
&& \int \psi(z) dL_{B_{n} +P_n}(z)  
   = \frac{1}{4\pi n} \int_{\C} \Delta \psi(z) \log \det(zI-B_n-P_n)(zI-B_n-P_n)^{*}
  dm(z)\\
  &&
  =\frac{1}{4\pi n} \intt_{\C} \Delta \psi(z) \log \det(|zI-B_n|-P_nU)(|zI-B_n|-P_nU)^{*}
  dm(z)
\end{eqnarray*}
where we used the
polar decomposition of $zI_n-B_n$ to write
$zI-B_n=|zI-B_n| U^*$ 
with $U$ a unitary matrix.  Since $P_nU$ has the same law 
as $P_n$, we are
back at the same setting as
in the proof of Theorem 
\ref{main-theo},
 with $|zI-B_n|$ replacing $T_n$.
It is then straightforward to check that the same 
arguments work under our present hypotheses;
the symmetrized  empirical  measure $\nu_n^z$  of  the singular values of 
$T_n(z)+P_n$
converges to $\tilde \Theta_z\boxplus\lambda_1$ by Lemma \ref{lem-S-conv}, 
which guarantees the convergence
of 
$$\int_\eps^{+\infty}\log |x|d\nu^z_n(x),$$ 
whereas our hypotheses allow us to bound uniformly the Stieltjes transform of 
$\nu^n_z$ on $\{z_1: \Im(z_1)\ge n^{-C_7}\}$ as in Lemma \ref{lem-tildegest}, 
hence providing a control of the integral on the interval $[n^{-C_7},\eps]$.
The control of the integral for $x<n^{-C_7}$ uses
a regularization by the Gaussian matrix $n^{-\gamma} N_n$
as in Proposition \ref{cor-D} .\qed

{\bf Acknowledgments:} We thank Greg Anderson for many
fruitful and encouraging discussions.
 We thank Yan Fyodorov 
 for pointing out the paper \cite{horn} and
  Philippe Biane
 for suggesting that our technique could be 
 applied to the examples in \cite{bl01}.
 We thank the referee for a careful reading of the manuscript.

\bibliographystyle{amsplain}

\begin{thebibliography}{10}
\bibitem{AGZ} Anderson, G. W., Guionnet, A. and Zeitouni, O.,
\textit{An introduction to random matrices}, Cambridge University Press,
Cambridge (2010). 
\bibitem{bai} Bai, Z., 
\textit{Circular law}, {Ann. Probab.} {\bf 25}, 494--529, (1997).
\bibitem{bl01}Biane, P. and Lehner, F., \textit{Computation of some examples of {B}rown's spectral measure in
              free probability} ,  {Colloq. Math.} {\bf 90}, 181--211 ,(2001). 
\bibitem{brown} Brown, L. G., \textit{Lidskii's theorem in the type II
case}, in ``Proceedings U.S.--Japan, 
Kyoto/Japan 1983'', Pitman Res. Notes. Math Ser. {\bf 123}, 1--35,
(1983).
\bibitem{CGM} Collins, B., Guionnet, A. and Maurel-Segala E.
\textit{Asymptotics of unitary and orthogonal  matrix integrals},
Adv. Math. {\bf 222}, 172--215, (2009).
\bibitem{FZ}
Feinberg, J. and
Zee, A.,
\textit{Non-Gaussian non-Hermitian random matrix theory: phase transition and
 addition formalism},
 Nuclear Phys. B {\bf   501}, 643--669,  (1997).

\bibitem{FS}
Fyodorov, Y.V. and Sommers, H.J.,
\textit{Spectra of random contractions and scattering theory for 
discrete-time systems},
{JETP Lett.} {\bf  72}, 422--426, (2000)

\bibitem{FW}
Fyodorov, Y.V. and Wei, Y.,
\textit{On the mean density of complex eigenvalues for an ensemble 
of random matrices with prescribed singular values. }
{Phys. A.} {\bf 41}, 502001,(2008).

\bibitem{ginibre} Ginibre, J.,  \textit{Statistical ensembles of complex, 
quaternion, and real matrices},
Jour. Math. Phys
{\bf 6}, 440--449, (1965).
\bibitem{girko}
    Girko, V. L.,
     \textit{The circular law},
   {Teor. Veroyatnost. i Primenen.}
    {\bf 29},
     {669--679},
      (1984).
\bibitem{gotzetikhomirov} G\"{o}tze, F. and Tikhomirov, A.,
\textit{The circular law for random matrices}, 	Annals Probab.
{\bf 38}, 1444--1491,
(2010).
\bibitem{Guionnet} Guionnet. A,
	\textit{Large random matrices: lectures on macroscopic asymptotics} 
	Lecture Notes in Mathematics {\bf 1957} 
	Lectures from the 36th Probability Summer School held in
              Saint-Flour, 2006, Springer-Verlag.
\bibitem{haageruplarsen} Haagerup, U.
and Larsen, F., 
\textit{Brown's spectral distribution measure for {$R$}-diagonal
 elements in finite von {N}eumann algebras}, 
J. Funct. Anal. {\bf 2}, 331--367, (2000).

\bibitem{HS} Haagerup, U.
and Schultz, A.,
\textit{Invariant subspaces for operators in a general
$II_1$ factor},
Publ. Math. IHES {\bf 109}, 19--111, (2009).
\bibitem{HS2} Haagerup, U.
and Schultz, A.,
\textit{Brown measures of unbounded operators 
affiliated with a finite von Neumann algebra},
Math. Scand. {\bf 100}, 209--263, (2007).

\bibitem{HT} Haagerup, U. and Thorbj{\o}rnsen, S. 
\textit{A new application of random matrices: {${\rm Ext}(C\sp *\sb
              {\rm red}(F\sb 2))$} is not a group},
  {Ann. of Math. (2)}
  {\bf 162}, 711--775, (2005).
    		



\bibitem{horn} Horn, A., 
\textit{On the eigenvalues of a matrix with prescribed singular
              values}, {Proc. Amer. Math. Soc.}, {\bf 5}, 4--7, (1954).

\bibitem{peresalsbook}
 Hough, J.B., Krishnapur, M., Peres, Y.  and Vir{\'a}g, B.,
\newblock {\em Zeros of Gaussian Analytic Functions and Determinantal Point
  Processes},
\newblock  Providence, RI, American Mathematical Society, (2009).

\bibitem{sommerslehman}
Lehmann, N.  and  Sommers, H.-J.,
\textit{Eigenvalue statistics of random real matrices},
Phys. Rev. Lett. {\bf 67}, 941--944, (2001).
\bibitem{NS} Nica, A. and Speicher, R.,
 \textit{${\mathcal R}$-diagonal pairs -- a common approach to 
Haar unitaries and circular elements},
 Fields Inst. Commun. {\bf 12}, 149--188 (1997).

      \bibitem{panzhou}
Pan, G. and Zhou,  W.,
\textit{Circular law, extreme singular values and potential theory},
J. Multivariate Anal. {\bf 101}, 645--656, (2010).
\bibitem{PS} Pastur, L. and Shcherbina, M.,
	\textit{Bulk universality and related properties of Hermitian
	matrix models}, J. Stat. Phys. {\bf 130}, 205--250, (2008).
\bibitem{SST} Sankar, A., Spielman, D. A. and Teng, S.-H., 
	\textit{Smoothed analysis of the conditioning number and
	growth factor of matrices}, 
	SIAM J. Matrix Anal. {\bf 28}, 446--476, (2006).
\bibitem{taovu} Tao, T. and Vu, V., 
\textit{Random matrices: the circular law},
 Commun. Contemp. Math.  {\bf 10},   261--307, (2008).
\bibitem{taovuk} Tao, T. and Vu, V., with appendix by M. Krishnapur,
\textit{Random matrices: Universality of ESD's and the circular law},
 	arXiv:0807.4898v5 [math.PR] (2008).
\bibitem{Vo91}
Voiculescu, D.,
\textit{Limit laws for random matrices and free products}
Inventiones Mathematicae {\bf 104}, 201--220, (1991).
\bibitem{VoiculescuAdvances}
Voiculescu, D.,
\textit{The analogues of entropy and of Fisher's information measure in free
 probability theory. VI. Liberation and mutual free information},
 Adv. Math.  {\bf 146}, 101--166, (1999).


\bibitem{weyl} Weyl, H.,
\textit{Inequalities between the two kinds of eigenvalues of a 
linear
transformation},  {Proc. Nat. Acad. Sci. U. S. A.} {\bf 35}, 408--411, (1949).
\bibitem{zysom} Zyczkowski, K. and  Sommers, H.-J.,
\textit{Truncations of random unitary matrices},
J. phys. A: Math. Gen. {\bf 33}, 2045--2057, (2000).

\end{thebibliography}

\end{document}